\documentclass[11pt,a4paper,reqno]{amsart}
\usepackage[T1]{fontenc}
\usepackage{amsmath}
\usepackage{amsthm}
\usepackage{amsfonts}
\usepackage{amssymb}
\usepackage{graphicx}
\usepackage{amsbsy}
\usepackage{mathrsfs}
\usepackage{bbm}
\newcommand{\be}{\begin{equation}}
\newcommand{\ee}{\end{equation}}
\newcommand{\ba}{\begin{array}}
\newcommand{\ea}{\end{array}}
\newcommand{\bea}{\begin{eqnarray}}
\newcommand{\eea}{\end{eqnarray}}
\newcommand{\bee}{\begin{eqnarray*}}
\newcommand{\eee}{\end{eqnarray*}}

\newtheorem{Thm}{Theorem}[section]
\newtheorem{Lemma}[Thm]{Lemma}
\newtheorem{Prop}[Thm]{Proposition}
\newtheorem{Definition}[Thm]{Definition}

\newtheorem{remark}[Thm]{Remark}
\numberwithin{equation}{section}

\catcode`@=11
\def\section{\@startsection{section}{1}%
  \z@{1.5\linespacing\@plus\linespacing}{.5\linespacing}%
  {\normalfont\bfseries\large\centering}}
\catcode`@=12

\def\C{{\mathbb C}}
\def\RR{\mathbb{R}}
\def\R{{\mathbb R}}

\def\lim{\mathop{\rm lim}}
\def\goto{\rightarrow}

\def\sup{\mathop{\rm sup}}

\def\e{\varepsilon}
\def\l{\lambda}

\def\log{{\rm log}}

\def\lsl{\frac{\lambda_s}{\lambda}}
\def\xsl{\frac{x_s}{\lambda}}
\def\tgamma{{\tilde{\gamma}}}

\def\qbnt{\tilde{Q}_b^{(0)}}

\def\er{\e_1}

\def\normeldloc{\int|\e|^2e^{-|y|}}
\def\normeldloch{\int|\eh|^2e^{-|y|}}

\def\bl{\mathcal L}

\def\qb{Q_b}

\def\pbnt{\tilde{P}_b^{(0)}}

\def\ct{{\mathcal{C}}^3}

\def\pbn{P_b^{(0)}}

\def\S{\Sigma}
\def\T{\Theta}

\def\dsdb{\frac{\partial \S}{\partial b}}
\def\dtdb{\frac{\partial \T}{\partial b}}
\def\dsdbh{\frac{\partial \Sh}{\partial b}}
\def\dtdbh{\frac{\partial \Th}{\partial b}}

\def\gab{\Gamma_b}

\def\zt{\hat{\zeta}}

\def\qbh{\hat{Q}_b}
\def\zh{\hat{\zeta}_b}
\def\zhre{\hat{\zeta}_1}
\def\zhim{\hat{\zeta}_2}
\def\Sh{\hat{\Sigma}}
\def\Th{\hat{\Theta}}

\def\eh{\hat{\e}}

\def\ut{\tilde{u}}

\def\Sh{\hat{\Sigma}}
\def\Th{\hat{\Theta}}

\def\fref#1{{\rm (\ref{#1})}}

\def\ds{\displaystyle}
\def\ni{\noindent}
\def\bs{\bigskip}

\title[]{Stable self similar blow up dynamics for slightly $L^2$ supercritical NLS equations}
\author[F. Merle]{Frank Merle}
\address{Universit\'e de Cergy Pontoise and IHES, France}
\email{frank.merle@math.u-cergy.fr}
\author[P. Rapha\"el]{Pierre Rapha\"el}
\address{IMT, Universit\'e Paul Sabatier, Toulouse, France}
\email{pierre.raphael@math.univ-toulouse.fr}
\author[J. Szeftel]{Jeremie Szeftel}
\address{CNRS, France, and Princeton University, USA}
\email{jszeftel@math.princeton.edu}

\begin{document}

\begin{abstract} We consider the focusing nonlinear Schr\"odinger equations $i\partial_t u+\Delta u +u|u|^{p-1}=0$ in dimension $1\leq N\leq 5$ and for slightly $L^2$ supercritical nonlinearities $p_c<p<(1+\e)p_c$ with $p_c=1+\frac{4}{N}$ and $0<\e\ll 1$. We prove the existence and stability in the energy space $H^1$ of a self similar finite time blow up dynamics and provide a qualitative description of the singularity formation near the blow up time
\end{abstract}

\maketitle


\section{Introduction}



\subsection{Setting of the problem}


We consider in this paper the nonlinear Schr\"odinger equation
\be
\label{nls}
(NLS) \ \  \left   \{ \begin{array}{ll}
         iu_t=-\Delta u-|u|^{p-1}u, \ \ (t,x)\in [0,T)\times \R^N\\
         u(0,x)=u_0(x), \ \ u_0:\R^N\to \C
         \end{array}
\right .
\ee
in dimension $1\leq N\leq 5$ with $$1<p<+\infty\textrm{ for }N=1, 2\textrm{ and }1<p<\frac{N+2}{N-2}\textrm{ for }N\geq 3.$$ From a result of Ginibre and Velo \cite{GV}, (\ref{nls}) is locally well-posed in $H^1=H^1(\R^N)$ and thus, for $u_0\in H^1$, there exists $0<T\leq +\infty$ and a unique solution $u(t)\in {\mathcal{C}}([0,T),H^1)$ to (\ref{nls}) and either $T=+\infty$, we say the solution is global, or $T<+\infty$ and then $\lim_{t\uparrow T}|\nabla u(t)|_{L^2}=+\infty$, we say the solution blows up in finite time.\\
(\ref{nls}) admits the following conservation laws in the energy space $H^1$: 
$$
\left   . \begin{array}{lll}
         L^2-\mbox{norm}: \ \ \int|u(t,x)|^2dx=\int|u_0(x)|^2dx;\\
         \mbox{Energy}:\ \ E(u(t,x))=\frac{1}{2}\int|\nabla u(t,x)|^2dx-\frac{1}{p+1}\int |u(t,x)|^{p+1}dx=E(u_0).\\
         \mbox{Momentum}:\ \ Im(\int\nabla u(t,x)\overline{u(t,x)}dx)=Im(\int\nabla u_0(x)\overline{u_0(x)}dx)
         \end{array}
\right .
$$
The scaling symmetry $\lambda^{\frac{2}{p-1}}u(\lambda^2t,\lambda x)$ leaves the homogeneous Sobolev space $\dot{H}^{\sigma_c}$ invariant with 
\be
\label{defsc}
\sigma_c=\frac{N}{2}-\frac{2}{p-1}.
\ee
From the conservation of the energy and the $L^2$ norm, the equation is subcritical for $\sigma_c<0$ and all $H^1$ solutions are global and bounded in $H^1$. The smallest power for which blow up may occur is $$p_c=1+\frac{4}{N}$$ which corresponds to $\sigma_c=0$ and is referred to as the $L^2$ critical case. The case $0<\sigma_c< 1$ is the $L^2$ super critical and $H^1$ subcritical case.\\

Even though the existence of finite time blow up dynamics for $\sigma_c\geq 0$ has been known since the 60' using standard global obstructive arguments based on the virial identity, \cite{ZS}, the explicit description of the singularity formation and of the different possible regimes is mostly open.


\subsection{The $L^2$ critical case and the stable log-log blow up }


 In the $L^2$ critical case $p=p_c$, a specific blow up regime of "log-log" type has been exhibited by Perelman \cite{P} in dimension $N=1$ and further extensively studied by Merle and Raphael in the series of papers \cite{MR1}, \cite{MR2}, \cite{R1}, \cite{MR3}, \cite{MR4}, \cite{MR5} where a complete description of this stable blow up dynamics is given together with sharp classification results in dimension $N\leq 5$.  The ground solitary wave solution to \fref{nls} which is the unique nonnegative radially symmetric solution to 
\be
\label{eqq}
\Delta Q_p-Q_p+Q_p^p=0, \ \ Q_p\in H^1,
\ee see \cite{GN}, \cite{KW}, plays a distinguished role in the analysis as it provides the blow up profile of log-log solutions.\\
Let $\Lambda$ be the generator of the scaling symmetry given by $\Lambda f=\frac{2}{p-1}f+y\nabla f$
 and recall the following Spectral Property which was proved in \cite{MR1} in dimension $N=1$ and checked numerically in dimensions $N\leq 5$ in \cite{FMR}:\\

{\bf Spectral Property} {\it Let $N\leq 5$ and $p=p_c$. Consider the two real Schr\"odinger operators \be\label{bln} \bl_1=-\Delta +\frac{2}{N}\left (\frac{4}{N}+1\right )Q_{p_c}^{\frac{4}{N}-1}y\cdot \nabla Q_{p_c}\ \ , \ \ \bl_2=-\Delta +\frac{2}{N}Q_{p_c}^{\frac{4}{N}-1}y\cdot \nabla Q_{p_c},\ee and the real valued quadratic form for $\e=\e_1+i\e_2\in H^1$: 
\be
\label{formquadrah}
H(\e,\e)=(\bl_1\e_1,\e_1)+(\bl_2\e_2,\e_2).
\ee 
Then there exists a universal constant $\delta_1>0$ such that $\forall \e \in H^1$, if $(\e_1,Q_{p_c})=(\e_1,\Lambda Q_{p_c})=(\e_1,yQ_{p_c})=(\e_2,\Lambda Q_{p_c})=(\e_2,\Lambda^2 Q_{p_c})=(\e_2,\nabla Q_{p_c})=0$, then $$H(\e,\e)\geq \delta_1 \left(\int |\nabla \e|^2+\int|\e|^2e^{-|y|}\right).$$}
We then have the following:

\begin{Thm}[Existence of a stable log-log regime, \cite{MR1}, \cite{MR2}, \cite{MR3}, \cite{R1}, \cite{MR4}, \cite{MR5}, \cite{FMR}]
\label{thmdeux}
Let $N\leq 5$ and $p=p_c$. There exists a universal constant $\alpha^*>0$ such that the following holds true. For any initial data $u_0\in H^1$ with small super-critical mass
\be
\label{estldeux}
|Q_p|_{L^2}<|u_0|_{L^2}<|Q_p|_{L^2}+\alpha^*
\ee
and nonpositive Hamiltonian $E (u_0)<0$,  the corresponding solution to (\ref{nls}) blows up in finite time $0<T<+\infty$ according to the following blowup dynamics: there exist geometrical parameters $(\lambda(t),x(t), \gamma(t))\in\RR^*_+\times\R^N\times\R$ and an asymptotic residual profile $u^*\in L^2$ such that: $$u(t)-\frac{1}{\lambda^{\frac{N}{2}}(t)}Q_p\left(\frac{x-x(t)}{\lambda(t)}\right)e^{i\gamma(t)}\to u^* \ \ \mbox{in} \ \ L^2.$$ The blowup point converges at blowup time: $$x(t)\to x(T) \in \R^N \ \ \mbox{as} \ \ t\to T,$$ the blowup speed is given by the log-log law
\be
\label{logloglawintro}\lambda(t)\sqrt{\frac{\log|\log(T-t)|}{T-t}}\to \sqrt{2\pi} \ \ \mbox{as} \ \ t\to T,\ee
and the residual profile satisfies:
 $$ u^*\in L^2 \ \ \mbox{but} \ \ u^*\notin L^p, \ \ \forall p>2.
 $$
 More generally, the set of initial data satisfying (\ref{estldeux}) and such that the corresponding solution to (\ref{nls}) blows up in finite time with the log-log law (\ref{logloglawintro}) is {\it open in $H^1$}.
\end{Thm}

In other words, the stable log-log regime corresponds to an {\it almost self similar regime} where the blow up solution splits into a singular part with a universal blow up speed and a universal blow up profile given by the exact ground state $Q_{p_c}$, and a regular part which {\it remains in the critical space} $u^*\in L^2$ but looses any regularity above scaling $u^*\notin H^{\sigma}$ for $\sigma>0$.


\subsection{On the super critical problem}


The explicit description of blow up dynamics in the super critical setting is mostly open. In fact, the only rigorous description of a blow up dynamics in a super critical setting is for $p=5$ in {\it any} dimension $N\geq 2$, see Raphael \cite{R2} for $N=2$,  Raphael, Szeftel \cite{RS} for $N\geq 3$. Note that this includes energy super critical problems. In this setting, the existence and {\it radial} stability of self similar solutions blowing up on an asymptotic blow up sphere -and not a blow up point- is proved. These solutions reproduce on the blow up sphere the one dimensional quintic and hence $L^2$ critical blow up dynamic and the blow up speed is indeed given by the log log law \fref{logloglawintro}.\\
In the formal and numerical work \cite{FGW}, Fibich, Gavish and Wang propose in dimension $N=2,3$ and for $p_c<p<5$ a generalization of the standing ring blow up solutions and investigate a blow up dynamic where the solution concentrates on spheres with radius collapsing to zero. Such dynamics are clearly exhibited numerically and seem to be stable by radial perturbations. A striking feature moreover is that these solutions {\it form a Dirac mass in $L^2$} like in the $L^2$ critical case: 
\be
\label{cnioeoe}
|u|^2\rightharpoonup M\delta_{x=0}+|u^*|^2
\ee for some universal quantum of mass $M>0$, and with a specific predicted blow up speed: 
\be
\label{blowupring}
|\nabla u(t)|_{L^2}\sim \frac{1}{(T-t)^{\frac{1}{1+\alpha}}}, \ \ \alpha=\frac{5-p}{(N-1)(p-1)}.
\ee
The rigorous derivation of such collapsing ring solutions is an important open problem in the field.\\

While ring solutions display a stability with respect to radial perturbations, they are widely believed to be unstable by non radial perturbations, \cite{FGW}. In fact, it has long been conjectured according to numerical simulations, see \cite{SS} and references therein, that the generic blow up dynamics in the super critical setting -at least for $p$ near $p_c$- should be of self similar type: $$u(t,x)\sim \frac{1}{[\lambda(t)]^{\frac{2}{p-1}}}P\left(\frac{x-x_T}{\lambda(t)}\right)e^{i\gamma(t)}$$ for some blow up point  $x_T\in \RR^N$ and where $\lambda$ is in the self similar regime: $$\lambda(t)\sim \sqrt{T-t}.$$ The delicate issue here is the profile $P$ which does not seem to be given by the ground state solution $Q_p$ to \fref{eqq} anymore. In fact, explicit self similar solutions may be computed using the following standard procedure. Let explicitly $$u(t,x)=\frac{1}{[\lambda(t)]^{\frac{2}{p-1}}}Q_b\left(\frac{x-x_T}{\lambda(t)}\right)e^{-\frac{\log(T-t)}{2b}}, \ \ b>0, \ \ \lambda(t)=\sqrt{2b(T-t)},$$ then $u$ solves \fref{nls} if and only if $Q_{b,p}=Q_b$ satisfies the nonlinear {\it stationary} elliptic PDE:
\be
\label{eppliticpf}
\Delta Q_b-Q_b+ib\Lambda Q_b+Q_b|Q_b|^{p-1}=0.
\ee
Exact {\it zero energy} solutions to \fref{eppliticpf} have been exhibited by Koppel and Landman, \cite{landman}, for slightly super critical exponent $p_c<p<(1+\e)p_c$, $\e\ll 1$, using geometrical ODE techniques. The existence of such solutions is related to a {\it nonlinear eigenvalue problem} and for all $p\in(p_c,(1+\e)p_c)$, a unique value $b(p)$ is found such that \fref{eppliticpf} admits a zero energy solution. An asymptotic law is derived which confirms previous formal computations, see \cite{SS} p147 and references therein: 
\be
\label{sigmac}
\sigma_c=e^{-\frac{\pi}{b(p)}(1+o(1))} \ \mbox{as} \ \ p\to p_c.
\ee
Moreover, the self similar profile converges to $Q_{p_c}$ locally: $$Q_{b(p)}\to Q_{p_c} \ \ \ \mbox{in} \ \ H^1_{loc} \ \ \mbox{as}  \ \ p\to p_c.$$ However the self similar  solutions belong to $\dot{H}^1\cap L^{p+1}$ but always miss the critical Sobolev space due to a logarithmic divergence at infinity: $$|Q_{b(p)}|(y)\sim \frac{1}{|y|^{\frac{N}{2}-\sigma_c}} \ \ \mbox{and hence} \ \ Q_{b(p)}\notin \dot{H}^{\sigma_c},$$ and thus they in particular miss $L^2$ and hence the physically relevant space $H^1$. Eventually, the construction of the self similar solution is delicate enough that it is not clear at all how this object should generate a {\it stable} self similar blow up dynamics for the time dependent problem


\subsection{Statement of the result}



The law for the nonlinear eigenvalue \fref{sigmac} is deeply related to the log-log law \fref{logloglawintro} which can be rewritten in the following form: $$\frac{ds}{ds}=\frac{1}{\lambda^2}, \ \ b=-\lsl, \ \ b_s= -e^{-\frac{\pi}{b}(1+o(1))} \ \ \mbox{as}\ \ s\to +\infty.$$ This intimate connection between the log-log law and the nonlinear eigenvalue problem underlying the self similar equation is at the heart of formal heuristics which first predicted \fref{sigmac}, we refer to the monograph \cite{SS} for a complete introduction to the history of the problem.\\

Our main claim in this paper is that the log-log analysis \cite{MR4} which allowed Merle and Raphael to derive the sharp log-log law for a large class of initial data provides a framework to somehow {\it bifurcate} from the critical value $p=p_c$ and prove the existence of a stable self similar regime in the energy space for slightly super critical exponents, with a blow up speed asymptotically satisfying the nonlinear eigenvalue relation \fref{sigmac}. In particular, our strategy  completely avoids the delicate issue of the existence of exact self similar solutions to \fref{eppliticpf}, and provides a framework to directly prove in a dynamical way the existence of {\it stable} self similar blow up dynamics, which we expect will apply to a large class of problems. In fact, we will show that it is enough to construct a {\it crude and compactly supported} approximation of the self similar profile \fref{eppliticpf}. Then, the intuition and technical tools inherited from the log-log analysis will allow us to obtain new {\it rigidity properties} and a {\it dynamical trapping} of the self similar regime. The outcome is a surprisingly robust proof of the existence of a $H^1$ stable self similar blow up regime.
\begin{Thm}[Existence and stability of a self similar blow up regime]
\label{thmmain}
Let $1\leq N\leq 5$. There exists $p^*>p_c$ such that for all $p\in(p_c,p^*)$, there exists $\delta(p)>0$ with $\delta(p)\to 0$ as $p\to p_c$, there exists $b^*(p)>0$ with 
\be
\label{loisigma}
\sigma_c=e^{-\frac{\pi}{b^*(p)}(1+\delta(p))}
\ee
and an open set $\mathcal O$ in $H^1$ of initial data such that the following holds true. Let $u_0\in \mathcal O$, then the corresponding solution to \fref{nls} blows up in finite time $0<T<+\infty$ according to the following dynamics: there exist geometrical parameters $(\lambda(t), x(t), \gamma(t))\in\RR^*_+\times\R^N\times\R$ and an excess of mass $\e(t)\in H^1$ such that: 
\be
\label{decomposition}
\forall t\in [0,T), \ \ u(t,x)=\frac{1}{\lambda^{\frac{2}{p-1}}(t)}[Q_p+\e(t)]\left(\frac{x-x(t)}{\lambda(t)}\right)e^{i\gamma(t)}
\ee
with
\be
\label{uniformbounde}
\ \ |\nabla \e(t)|_{L^2}\leq \delta(p).
\ee
The blowup point converges at blowup time: 
\be
\label{convblowuppiint}
x(t)\to x(T) \in \R^N \ \ \mbox{as} \ \ t\to T,
\ee 
and the blow up speed is self similar: 
\be
\label{selfsmilar}
\forall t\in [0,T), \ \ (1-\delta(p))\sqrt{2b^*(p)}\leq \frac{\lambda(t)}{\sqrt{T-t}}\leq (1+\delta(p))\sqrt{2b^*(p)}.
\ee 
Moreover, there holds the strong convergence:
\be
\label{convustarhsigma}
\forall \sigma\in[0,\sigma_c), \ \ u(t)\to u^* \ \ \mbox{in} \ \ H^{\sigma} \ \ \mbox{as}\ \ t\to T,
\ee 
and the asymptotic profile $u^*$ displays a singular behavior on the blow up point: $\exists R(u_0), C(p)>0$ such that:
 \be
 \label{nonsmoothustar}
\forall R\in (0,R(u_0)), \ \  C(p)(1-\delta(p))\leq \frac{1}{R^{2\sigma_c}}\int_{|x|\leq R}|u^*|^2\leq C(p)(1+\delta(p)).
\ee In particular: $$u^*\in H^{\sigma} \ \ \mbox{for} \ \ 0\leq \sigma<\sigma_c \ \ \mbox{but} \ \ u^*\notin \dot{H}^{\sigma_c}.$$

\end{Thm}

{\it Comments on Theorem \ref{thmmain}}\\

{\it 1. Asymptotic dynamics}: Let us stress onto the fact that \fref{decomposition} does not give a sharp asymptotic on the singularity formation, and in particular the question of the possible asymptotic stability in renormalized variables in open. In fact, we will construct a rough approximate self similar profile $b\to Q_b$ and show that the solution decomposes into $$u(t,x)=\frac{1}{\lambda^{\frac{2}{p-1}}(t)}[Q_{b(t)}+\e]\left(\frac{x-x(t)}{\lambda(t)}\right)e^{i\gamma(t)}$$ with $-\lambda_t\lambda\sim b(t)$ and $|\e(t)|_{\dot{H}^1}\ll 1$. The key will be to prove a dynamical trapping on the projection parameter $b(t)$ and a uniform control on the radiation $\e$: $$\forall t\in [0,T),  \ \ (1-\delta(p))b^*(p)\leq b(t)\leq (1+\delta(p))b^*(p) \ \ \mbox{and} \ \ |\e(t)|_{\dot{H}^1}\ll 1 ,$$  but this does not exclude possible oscillations of both $b(t)$ and $\e(t)$. In some sense, this confirms the analysis in \cite{MR1}, \cite{MR2}, see also Rodnianski, Sterbenz \cite{RodSter}, Raphael, Rodnianski \cite{RRS}, where the key observation was that is is not necessary to obtain a complete description of the dispersive structure of the problem in renormalized variables to prove finite time blow up. The main difference however with these works is that we are here in a situation where there holds no a priori orbital stability bound on $\e(t)$, and a spectacular feature is that rough profiles are enough to capture the self similar blow up speed. This somehow confirms the robustness of the log-log analysis developed in \cite{MR2}, \cite{MR3}, \cite{MR4}. Our perturbative approach and the strategy of bifurcating from the critical case is also reminiscent from the general framework developed by Fibich and Papanicolaou, \cite{FP}.\\

{\it 2. On the behavior of the critical norm}: In \cite{MR7}, Merle and Raphael showed in the range of parameters $N\geq 3$, $0<\sigma_c<1$, that any radially symmetric finite time blow up solution in $H^1$ must leave the critical space at blow up time with a lower bound: 
\be
\label{nknoe}
|u(t)|_{\dot{H}^{\sigma_c}}\geq |\log(T-t)|^{\alpha(p)}\ \ \mbox{as} \ \ t\to T.
\ee 
The self similar solutions constructed from Theorem \ref{thmmain} satisfy a logarithmic upper bound --see Remark \ref{rkintro}--:
\be
\label{cnnooeen}
|u(t)|_{\dot{H}^{\sigma_c}}\lesssim |\log(T-t)|^{\frac{p}{p+1}}\ \ \mbox{as} \ \ t\to T
\ee which proves the sharpness in the logarithmic scale of the lower bound \fref{nknoe}. A logarithmic lower bound could also be derived after some extra work. The singularity \fref{nonsmoothustar} in fact sharpens for this specific class of initial data the divergence \fref{cnnooeen}, and is according to \fref{nonsmoothustar} the major difference between the $L^2$ critical blow up where the $L^2$ conservation law forces the radiation to remain in the critical $L^2$ space, and the super critical blow up where the radiation leaves asymptotically $\dot{H}^{\sigma_c}$.\\

{\it 3. Weak blow up}: On the contrary to the ring solutions which concentrate in $L^2$ at blow up time according to \fref{cnioeoe}, the self similar solutions constructed from Theorem \ref{thmmain} do not concentrate mass according to \fref{convustarhsigma}, a situation which is referred to as {\it weak} blow up, see \cite{SS}. A striking feature of the analysis is that the $L^2$ conservation law, even though below scaling, plays a fundamental role in the proof of the stabilization of the self similar blow up. 


\subsection{Strategy of the proof}


Let us briefly state the main steps of the proof of Theorem \ref{thmmain}.\\

{\bf step 1}  Construction of an approximate self similar profile.\\

Let a small parameter $b>0$ and recall that the equation for self similar profiles \fref{eppliticpf} does not admit solutions in $H^1$. The first step is to construct an approximate solution $Q_b$ which is essentially compactly supported and satisfies an approximate self similar equation: 
$$\Delta Q_b-Q_b+ib\Lambda Q_b+Q_b|Q_b|^{p-1}=\Psi_b+O(\sigma_c^2),$$ see Proposition \ref{propapproximatesolution}. This profile incorporates the leading order $O(\sigma_c)$ deformation with respect to the $L^2$ critical profiles constructed in \cite{MR2}, \cite{MR4}, while the error $\Psi_b$ is a far away localized error which is inherited from the space localization of the profile to avoid the slowly decaying tails.\\

{\bf step 2} Dynamical trapping of $b$.\\

We now chose initial data such that on a short time, the solution admits a decomposition $$u(t,x)=\frac{1}{\lambda^{\frac{2}{p-1}}(t)}(Q_{b(t)}+\e)\left(t,\frac{x-x(t)}{\lambda(t)}\right)e^{i\gamma(t)}\ \  \mbox{with} \ \ |\nabla \e(t)|^2_{L^2}\ll e^{-\frac{\pi}{b(t)}}$$ and aim at deriving a dynamical trapping for the parameter $b(t)$ and the deformation $\e(t)$. More precisely, introducing the global rescaled time $\frac{ds}{dt}=\frac{1}{\lambda^2}$, the dynamical system driving $\lambda$ is $$-\lambda\lambda_t=-\lsl\sim b(t)$$ and hence finite time blow up in the self similar regime will follow from $$b(t)\sim b_0>0, \ \ |\nabla \e(t)|^2_{L^2}\ll e^{-\frac{\pi}{b(t)}(1-c)}$$ for some small constant $c>0$ in the maximum time interval of existence, see Proposition \ref{propboot}. In order to derive such dynamical controls, we run the log-log analysis developed in \cite{MR4} by keeping track of the leading order $O(\sigma_c)$ deformation. The outcome is the derivation of two structural monotonicity formulae for the parameter $b$. The first one is inherited from the local virial control first derived in \cite{MR1} and roughly leads to $$\sigma_c+|\nabla \e|_{L^2}^2-e^{-\frac{\pi}{b}}\lesssim b_s,$$ see Proposition \ref{propviriel}. The positive term $+\sigma_c$ in the above LHS is a non trivial supercritical effect and is a consequence of the structure of the self similar profiles in the local range $|y|\leq 1$. The second monotonicity formula is a consequence of the $L^2$ conservation law and the control of the mass ejection phenomenon: $$b_s\lesssim \sigma _c-e^{-\frac{\pi}{b}},$$ where the nonpositive term $-e^{-\frac{\pi}{b}}$ is obtained from a {\it flux computation in the far away zone $|y|>>1$} which is based on the presence of the slow decaying tails of self similar solutions. The outcome is the derivation of the dynamical system for $b$: $$b_s\sim \sigma_c-e^{-\frac{\pi}{b}}$$ which traps $b$ around the value $$b\sim b^* \ \ \mbox{with} \ \ \sigma_c\sim e^{-\frac{\pi}{b^*}}.$$ Note that this shows that the self similar blow up speed is derived from the constraints both on compact sets and at infinity where the dispersive mass ejection process is submitted to the global constraint of the $L^2$ conservation law. Note also that the $L^2$ critical log-log law corresponds to the dynamical system $$b_s\sim -e^{-\frac{\pi}{b}}$$ and hence the supercritical self similar blow up appears as directly branching from the $L^2$ critical $\sigma_c=0$ degenerate log-log blow up. The nontrivial {\it pointwise} control on $\e$ $$|\nabla \e(t)|^2_{L^2}\ll e^{-(1-c)\frac{\pi}{b(t)}}$$ will also follow from the obtained Lyapounov controls. Here a difficulty will occur with respect to the $L^2$ critical case to control the nonlinear terms in $\e$ due in particular to the unboundedness of the scaling invariant critical Sobolev norm. A new strategy inspired from  \cite{RS}, \cite{RodSter}, \cite{RRS} is derived  which relies on the control of Sobolev norms {\it strictly above scaling} in the self similar regime, see section \ref{sectionclosingboot}.\\
The conclusions of Theorem 2 are now a simple consequence of this dynamical trapping of the solution.\\

This paper is organized as follows. In section \ref{sectiontwo}, we construct approximate self similar solutions, Proposition \ref{propapproximatesolution}. We then describe the set of initial data leading to self similar blow up, Definition \ref{defdefp}, and set up the bootstrap argument, Proposition \ref{propboot}. In section \ref{sectionthree}, we derive the key dynamical controls and the two monotonicity formulae, Proposition \ref{propviriel} and Proposition \ref{lemmadiffineqb}. In section \ref{sectionproofthm}, we close the bootstrap argument as a consequence of the obtained Lyapounov type controls and conclude the proof of Theorem \ref{thmmain}.

\bs
\ni

\noindent{\bf Acknowledgments.} F.M. is supported by ANR Projet Blanc OndeNonLin. P.R and J.S are supported by ANR jeunes chercheurs SWAP. \\

\noindent{\bf Notations}  We let $|\nabla|^{\sigma}$ be the Fourier multiplier $\widehat{|\nabla|^{\sigma}f}(\xi)=|\xi|^{\sigma}\widehat{f}(\xi)$. We let $Q_p$ be the unique radially symmetric nonnegative solution in $H^1$ to $$\Delta Q_p-Q_p+Q_p^{p+1}=0.$$ We introduce the error to $L^2$ criticality:
\be
\label{defalpha}
p_c=1+\frac{4}{N}, \ \ \sigma_c=\frac{N}{2}-\frac{2}{p-1}=\frac{N(p-p_c)}{2(p-1)},
\ee
where $0<\sigma_c<1$ is the Sobolev critical exponent. We introduce the associated scaling generators:
\be
\label{genrtor}
\Lambda f=\frac{2}{p-1}f+y\cdot\nabla f, \ \ Df=\frac{N}{2}f+y\cdot\nabla f=\Lambda f+\sigma_c f.
\ee 
We denote the $L^2(\RR^N)$ scalar product $$(f,g)=\int_{\RR^N}f(x)g(x)dx$$ and observe the integration by parts formula: 
\be
\label{integraoarts}
(Df,g)=-(f,Dg), \ \ (\Lambda f,g)=-(f,\Lambda g+2\sigma_c g).
\ee
We let $L=(L_+,L_-)$ be the linearized operator close to $Q_p$: 
\be
\label{deflpluslmoins}
L_+=-\Delta +1-pQ_p^{p-1}, \ \ L_-=-\Delta +1-Q_p^{p-1}.
\ee


\section{Description of the blow up set of initial data}
\label{sectiontwo}


This section is devoted to the description of the open $H^1$ set $\mathcal O$ of initial data leading to the self similar blow up solutions described by Theorem \ref{thmmain} which relies on the construction of  approximate self similar profiles.


\subsection{Construction of approximate self similar solutions}


Our aim in this section is to construct suitable approximate solutions to the self similar equation \fref{eppliticpf}. Let us make the following general ansatz: $$ u(t,x)=\frac{1}{\lambda^{\frac{2}{p-1}}(t)}v\left(t,\frac{x-x(t)}{\lambda(t)}\right)e^{i\gamma(t)}$$ and introduce the rescaled time $$\frac{ds}{dt}=\frac{1}{\lambda^2(t)},$$ then $u$ is a solution to (\ref{nls}) if and only if $v$ solves: $$i\partial_s v+\Delta v-v-i\lsl\Lambda v +v|v|^{p-1}=(\gamma_s-1)v+i\frac{x_s}{\lambda}\cdot\nabla v.$$ Let us fix $$\gamma_s=1, \ \ x_s=0, \ \ -\lsl=b(s)$$ and look for solutions of the form: $$v(s,y)=Q_{b(s)}(y)$$ where the unknowns are the mappings $$s\to b(s), \ \ \ b\to Q_b.$$ This corresponds a to a {\it slow variable} formulation of a generalized self similar equation which goes back to previous formal works, see \cite{SS}, and was rigorously used in \cite{P}, \cite{MR2}, see also \cite{KMR}, \cite{RRS} for related transformations in different settings. To prepare the computation, we let $$b_s=\sigma_c\mu_b$$
where $\mu_b$ is a function of $b$ which will be made explicit later on. The generalized self similar equation becomes: 
\be
\label{selfsimilqb}
i\sigma_c \mu_b\frac{\partial Q_b}{\partial b}+\Delta Q_b-Q_b+ib\Lambda Q_b+Q_b|Q_b|^{p-1}=0.
\ee 
Let us perform the conformal change of variables $$P_b=Q_be^{\frac{ib|y|^2}{4}},$$ then a simple algebra leads to:
\be
\label{eqpbgenerale}
i\sigma_c\mu_b\frac{\partial P_b}{\partial b}+\Delta P_b-P_b-i\sigma_cb P_b+\frac{1}{4}(b^2+\sigma_c\mu_b)|y|^2P_b+P_b|P_b|^{p-1}=0.
\ee
Our aim is to {\it find $\mu_b$} so as to be able to construct an approximate solution to \fref{eqpbgenerale} with an error of order formally $\sigma_c^2$ in the region $|y|\leq \frac{1}{b}$. Our construction is elementary and relies on the computation of the first term in the Taylor expansion of $(\mu_b,Q_b)$ in $\sigma_c$ near the $L^2$ critical value $\sigma_c=0$.\\

Let a small parameter $0<\eta<<1$ to be fixed later, a non zero number $b$, and set
\be
\label{rplusrmoins}
 R_b=\frac{2}{|b|}\sqrt{1-\eta}\ \ \mbox{and} \ \ R_b^-=\sqrt{1-\eta} R_b .
\ee
Denote $B_{R_b}=\{y\in \R^N, \ \ |y|< R_b\}$ and $\partial B_{R_b}=\{y\in \R^N, \ \ |y|= R_b\}$. We introduce a regular radially symmetric cut-off function 
\be
\label{defcutphib}
\phi_b(x)=\left\{\begin{array}{ll} 0 \ \ \mbox{for} \ \ |x|\geq R_b,\\
						1 \ \ \mbox{for} \ \ |x|\leq R_b^-,
						\end{array}
						\right .
						\ \ 0\leq \phi_b(x)\leq 1,
						\ee 
						such that:
\be
\label{derivephib}
|\phi'_b|_{L^{\infty}}+|\Delta \phi_b|_{L^{\infty}}\to 0  \ \ \mbox{as} \ \ |b|\to 0.
\ee
We also consider the norm on radial functions $\|f\|_{{\mathcal{C}}^j}=\max_{0\leq k\leq j}\|f^{(k)}(r)\|_{L^{\infty}(\R_+)}$. Let us start with the construction of the profile for $\sigma_c=0$ where we view $\sigma_c$ and $p$ as independent parameters and leave $p$ supercritical.

\begin{Prop}[$Q_b^{(0)}$ profiles]
\label{profilespbzero}
There exists $p^*>p_c$ and $C,\eta^*>0$ such that for all $p_c\leq p<p^*$, for all $0<\eta<\eta^*$, there exists $b^*(\eta), \e^*(\eta)>0$ going to zero as $\eta\to 0$ such that for all $|b|\leq b^*(\eta)$, the following holds true:\\
(i) Existence of a unique $P_b$ profile: there exists a unique radial solution $P^{(0)}_b$ to 
\be
\label{eqpbzero}
\ \  \left   \{ \begin{array}{lll}
       \Delta P^{(0)}_b -\pbn+\frac{b^2|y|^2}{4}\pbn+\left(\pbn\right)^{p}=0,\\
P_b^{(0)}>0 \ \ \mbox{in} \ \ B_{R_b},\\
        P^{(0)}_b(0)\in (Q_p(0)-\e^*(\eta), Q_p(0)+\e^*(\eta)), \ \ P^{(0)}_b(R_b)=0.
         \end{array}
\right .
\ee
Moreover, let 
\be
\label{defpbtilde}
\pbnt(r)=\pbn(r)\phi_b(r),
\ee 
then $\pbnt$ is twice differentiable with respect to $b^2$ with uniform estimate:
\be
\label{convunifqb}
\left\|e^{(1-\eta)\frac{\theta(|b|r)}{|b|}}\left(\pbnt-Q_p\right)\right\|_{\ct}+\left\|e^{(1-\eta)\frac{\theta(|b|r)}{|b|}}\left(\frac{\partial \pbnt}{\partial b^2}-\rho\right)\right\|_{{\mathcal{C}}^2}\to 0 \ \ \mbox{as} \ \ b\to 0,
\ee
\be
\label{convunifqbbis}
\left\|e^{(1-\eta)\frac{\theta(|b|r)}{|b|}}\frac{\partial^2\tilde{P}_b^{(0)}}{\partial^2(b^2)}\right\|_{\ct}\leq C
\ee
where 
\be
\label{deftheta}
\theta(w)=1_{0\leq w\leq 2}\int_{0}^{w}\sqrt{1-\frac{z^2}{4}}dz+1_{w> 2}\frac{\theta(2)}{2}w,
\ee
and $\rho$ is the unique solution in $H^1_{rad}$ to 
\be
\label{defrho}
L_+\rho=\frac{1}{4}|y|^2Q.
\ee
(ii) Properties of the $\qbnt$ profile: $\qbnt=e^{-ib\frac{|y|^2}{4}}\pbnt$ satisfies:
\be
\label{eqqbfhji}
\Delta\qbnt -\qbnt+ibD\qbnt+\qbnt|\qbnt|^{p-1}=-\Psi^{(0)}_b=-\tilde{\Psi}_b^{(0)}e^{ib\frac{|y|^2}{4}}
\ee
\be
\label{psib}
-\tilde{\Psi}^{(0)}_b=2\nabla\phi_b\cdot\nabla\pbn+\pbn (\Delta \phi_b)+(\phi_b^{p}-\phi_b)(\pbn)^p
\ee
and for any polynomial $f(y)$ and integer $k=0,1$, 
\be
\label{unpreciseest}
\left|f(y)\frac{\partial \tilde{\Psi}_b^{(0)}}{\partial y^k}\right|_{L^{\infty}}\leq e^{-\frac{C_{p}}{|b|}}
\ee
for some constant $c_p$ depending on $p$.\\
(iii) Computation of the momentum and the $L^2$ mass: there holds
\be
\label{cnenveo}
Im\left(\int \nabla\qbnt\overline{\qbnt}\right)=0, \ \ Im\left(\int y\cdot\nabla\qbnt\overline{\qbnt}\right)=-\frac{b}{2}|y\qbnt|_2^2
\ee and the supercritical mass property:
\be
\label{qbtsupcrit}
\frac{d^2}{db^2}\left(\int|\qbnt|^2\right)_{|b=0}=c_0(p)\ \ \mbox{with} \ \ c_0(p)\to c_0(p_c)>0 \ \ \mbox{as} \ \ p\to p_c.
\ee
\end{Prop}

The proof of Proposition \ref{profilespbzero} is parallel to the one of Proposition in \cite{MR2} and Proposition \cite{MR3} and relies on standard elliptic techniques and the knowledge of the kernel of the linearized operator close to Q, explicitly: 
\be
\label{rappelkernelL}
Ker(L_+)=\mbox{span}(\nabla Q), \  \ Ker(L_-)=\mbox{span}(Q),
\ee 
see \cite{weinstein}, \cite{spectrenls}, and the fact that we are working here before the turning point $\frac{2}{b}$ and hence with uniformly elliptic operators. The detailed proof is left to the reader.

\begin{remark} The question of the value of the energy of the modified profile is an important issue. It will indeed be computed in Proposition \ref{propapproximatesolution} for the full approximate profile, see \fref{calculhamiltonian}.
\end{remark}

After the turning point $\frac{2}{b}$, leading order phenomenons are of linear type and a natural prolongation of the approximate blow up profile is given by the so called outgoing radiation. 

\begin{Lemma}[Linear outgoing radiation]
\label{outgoingradiation}
See Lemma 15 in \cite{MR3}. There exist universal constants $C>0$ and $\eta^{*}>0$ such that $\forall 0<\eta<\eta^*$, there exists $b^{*}(\eta)>0$ such that $\forall 0<b<b^{*}(\eta)$, the following holds true: let $\Psi^{(0)}_b$ be given by (\ref{eqqbfhji}), there exists a unique radial solution $\zeta_b$ to  
\be
\label{zeta}
\ \  \left   \{ \begin{array}{ll}
         \ds\Delta\zeta_b-\zeta_b+ibD\zeta_b=\Psi^{(0)}_b,\\
         \ds\int|\nabla \zeta_b|^2<+\infty .
         \end{array}
\right .
\ee
Moreover, let $\theta$ be given by (\ref{deftheta}), and consider 
\be
\label{defgab}
\gab=\lim_{|y|\to+\infty}|y|^N|\zeta_b(y)|^2,
\ee
then there holds: 
\be
\label{estgammainfinity}
\left| |y|^{\frac{N}{2}}(|\zeta_b|+|y||\nabla(\zeta_b)|)\right|_{L^{\infty}(|y|\geq R_b)}\leq \gab^{\frac{1}{2}-C\eta}<+\infty,\ee
\be
\label{estgradzetab}
\int |\nabla \zeta_b|^2\leq \gab^{1-C\eta}.
\ee 
For $|y|$ large, we have more precisely:
\be
\label{gammapositif}
\forall |y|\geq R_b^2, \ \ e^{-2(1-C\eta)\frac{\theta(2)}{b}}\geq |y|^{N}|\zeta_b(y)|^2\geq \frac{4}{5}\gab\geq e^{-2(1+C\eta)\frac{\theta(2)}{b}},
\ee
\be
\label{controlegrad}
\forall |y|\geq R_b^2, \ \ |\nabla \zeta_b(y)|\leq \frac{C}{|y|^{1+\frac{N}{2}}}\frac{\gab^{\frac{1}{2}}}{|b|}.
\ee
For $|y|$ small, we have: $\forall \sigma\in(0,5)$, $\exists \eta^{**}(\sigma)$ such that $\forall 0<\eta<\eta^{**}(\sigma)$, $\exists b^{**}(\eta)$ such that $\forall 0<b<b^{**}(\eta)$, there holds:
\be
\label{estgammacompact}
\left |\zeta_b(y)e^{-\sigma\frac{\theta(b|y|)}{b}}\right|_{{\mathcal{C}}^2(|y|\leq R_b)}\leq \gab^{\frac{1}{2}+\frac{1}{10}\sigma}.\ee
Last, $\zeta_b$ is differentiable with respect to $b$ with estimate
\be
\label{estdiffzetab}
\left|\frac{\partial\zeta_b}{\partial b}\right|_{{\mathcal{C}}^1}\leq \gab^{\frac{1}{2}-C\eta}.
\ee
\end{Lemma}

\begin{remark} Recall from \fref{deftheta} that $\theta(2)=\frac{\pi}{2}$. Moreover it is enough for our analysis to compute the radiation with the $L^2$ scaling generator $D$ in \fref{zeta} instead of $\Lambda$ as the error will generate lower order terms.
\end{remark}

\begin{remark}
Now that $\Gamma_b$ has been defined in \eqref{defgab}, we can give a more precise formulation of the estimate \eqref{unpreciseest} satisfied by $\tilde{\Psi}^{(0)}_b$ which is a direct consequence of \fref{convunifqb} and the localization procedure \fref{psib}: for any polynomial $f(y)$ and integer $k=0,1$, 
\be
\label{preciseest}
\left|f(y)\frac{\partial \tilde{\Psi}_b^{(0)}}{\partial y^k}\right|_{L^{\infty}}\leq \Gamma_b^{\frac{1}{2}(1-C\eta)}.
\ee
\end{remark}

The proof of Lemma \ref{outgoingradiation} is completely similar to the one of Lemma 15 in \cite{MR3} -see also Appendix A in \cite{MR4}- and hence left to the reader. In particular, the nondegeneracy \fref{gammapositif} is a consequence of the {\it nonlinear} construction of the profile $\pbnt$.\\

One should think of $\qbnt$ as being an approximate solution to the generalized self similar equation \fref{selfsimilqb} with $\mu_b=0$ -self similar law- and an error of order $\sigma_c$ in the elliptic zone $|y|\leq \frac{1}{b}$. We now claim that there exists a -locally unique- non trivial $\mu_b(p)>0$ which allows one to construct an approximate solution to the generalized self similar equation \fref{selfsimilqb} of order $O(\sigma_c^2+\gab^2)$ in the elliptic zone $|y|\leq \frac{2}{b}$:

\begin{Prop}[Approximate generalized self similar profiles]
\label{propapproximatesolution}
There exists $p^*>p_c$ and $C,c_1,\eta^*>0$ such that for all $p_c<p<p^*$ and $0<\eta<\eta^*$, there exists $c_0(p)>0$, $b^*(\eta)>0$ such that for all $|b|<b^*(\eta)$, the following holds true:\\
(i) Construction of the modified profile: there exist $\mu_b=\mu(b,p)>0$ and a radially symmetric complex valued function $T_b=T(b,p)$ with 
\be
\label{uniformestimatest}
\left\|e^{(1-\eta)\frac{\theta(|b|r)}{|b|}}T_b(r)\right\|_{\ct}+\left\|e^{(1-\eta)\frac{\theta(|b|r)}{|b|}}\frac{\partial T_b(r)}{\partial b}\right\|_{\mathcal C^2}+\left|\frac{\partial \mu_b}{\partial b}\right|\leq C \  \ \mbox{as} \ \ b\to 0,
\ee
\be
\label{estmubp}
\mu_b\to \frac{8|Q_p|_{L^2}^2}{(1+2\sigma_c)|yQ_p|_{L^2}^2}  \ \ \mbox{as} \ \ b\to 0
\ee
with the following properties. Let $$P_b=\tilde{P}_b^{(0)}+\sigma_c T_b, \ \ Q_b=P_be^{i\frac{b|y|^2}{4}},$$ then 
\be
\label{defpsib}
\Psi_b=-i\sigma_c \mu_b\frac{\partial Q_b}{\partial b}-\Delta Q_b+Q_b-ib\Lambda Q_b-Q_b|Q_b|^{p-1}=\Psi_b^{(0)}+\Psi_b^{(1)}
\ee
with $\Psi_b^{(0)}$ given by \fref{eqqbfhji} satisfies \eqref{preciseest}, and for $k=0,1$:
\be
\label{etglobalepsitwo}
\left\|e^{(1-\eta)\frac{\theta(|b|r)}{|b|}}\frac{\partial^k \Psi_b^{(1)}}{\partial y^k}\right\|_{L^{\infty}(|y|\leq R_b^-)}\leq \sigma_c^{1+c_1},
\ee
for some constant $c_1>0$ depending only on $N$ and:
\be
\label{etglobalepsitwobis}
\left\|e^{(1-\eta)\frac{\theta(|b|r)}{|b|}}\frac{\partial^k \Psi_b^{(1)}}{\partial y^k}\right\|_{L^{\infty}(|y|\geq R_b^-)}\leq C\sigma_c. 
\ee
(ii) Estimate of the invariants on $Q_b$: there holds
\be
\label{degenmoment} 
Im\left(\int \nabla Q_b\overline Q_b\right)=0,
\ee
\be
\label{degenviriel}
Im\left(\int y\cdot\nabla Q_b\overline Q_b\right)=-\frac{b}{2}|yQ_p|_2^2(1+O(|b|+\sigma_c)) \ \ \mbox{as} \ \ b \to 0,
\ee
\be
\label{qbtsupcritbis}
\int |Q_b|^2=\int|Q_p|^2+M(b)+O(\sigma_c)
\ee
with 
\be
\label{propmb}
M(0)=0 \ \ \mbox{and} \ \  \frac{d^2}{db^2}M(b)|_{b=0}=c_0(p)\to c_0(p_c)>0 \ \ \mbox{as} \ \ p\to p_c,
\ee
and the degeneracy of the Hamiltonian:
\be
\label{calculhamiltonian}
\left|E(Q_b)\right|\leq \gab^{1-C\eta}+C\sigma_c.
\ee

\end{Prop}

{\bf Proof of Proposition \ref{propapproximatesolution}}\\

The proof relies on a Taylor expansion in $\sigma_c$ of formal solutions to \fref{selfsimilqb}. The choice of $\mu_b$ is dictated by the presence of a non trivial kernel for the operator $L_-$ driving the imaginary part of \fref{selfsimilqb} near $Q_p$ as given by \fref{rappelkernelL}, while $L_+$ on the real part is invertible in the radial sector. The computation of the invariants and in particular the energy degeneracy \fref{calculhamiltonian} follow from Pohozaev identity.\\

{\bf step 1} Construction of $T_b$.\\

Let $Q_b=P_be^{-ib\frac{|y|^2}{4}}$ and $\Psi_b$ given by \fref{defpsib} , then: $\Psi_b=\tilde{\Psi}_be^{-ib\frac{|y|^2}{4}}$ with:
\be
\label{defpsitilde}
-\tilde{\Psi}_b=i\sigma_c\mu_b\frac{\partial P_b}{\partial b}+\Delta P_b-P_b-i\sigma_c b P_b+\frac{1}{4}(b^2+\sigma_c\mu_b)|y|^2P_b+P_b|P_b|^{p-1}.
\ee
Let $\tilde{\Psi}^{(0)}_b$ be given by \fref{psib}, equivalently:
\be
\label{psitildepbo}
-\tilde{\Psi}^{(0)}_b=\Delta \pbnt-\pbnt+\frac{b^2}{4}|y|^2\pbnt+\pbnt|\pbnt|^{p-1}.
\ee
We expand $P_b=\tilde{P}_b^{(0)}+\sigma_c T_b$. Let $\phi_b$ be the cut off function given by \fref{defcutphib}. We compute:
\bea
\label{caclulpsitilde}
\nonumber -\tilde{\Psi}_b& = & -\tilde{\Psi}_b^{(0)}+i\sigma_c^2\mu_b\frac{\partial T_b}{\partial b}-i\sigma_c^2 bT_b+\frac{1}{4}\sigma_c^2\mu_b|y|^2T_b+\frac{1}{4}\sigma_c b^2(1-\phi_b)|y|^2T_b\\
\nonumber& + & (\pbnt+\sigma_c T_b)|\pbnt+\sigma_c T_b|^{p-1}-(\pbnt)^{p}-\sigma_c p(\pbnt)^{p-1}Re(T_b)-i\sigma_c (\pbnt)^{p-1}Im(T_b)\\
& + & \sigma_c\left[-(L_{b})_+Re(T_b)+\frac{\mu_b}{4}|y|^2\pbnt\right]+ i\sigma_c\left[-(L_{b})_-Im(T_b)+\mu_b\frac{\partial \pbnt}{\partial b}-b\pbnt\right]
\eea
where we introduced the linearized operators close to $\pbnt$:
$$(L_+)_b=-\Delta +1-\frac{b^2}{4}\phi_b|y|^2-p(\pbnt)^{p-1}, \ \ (L_-)_b=-\Delta +1-\frac{b^2}{4}\phi_b|y|^2-(\pbnt)^{p-1}.$$ We thus aim at finding $(\mu_b,T_b)$ so as to cancel the $O(\sigma_c)$ in the RHS of \fref{caclulpsitilde}: 
\be
\label{eqlpluslmoins}
-(L_{b})_+Re(T_b)+\frac{\mu_b}{4}|y|^2\pbnt=0, \ \ -(L_{b})_-Im(T_b)+\mu_b\frac{\partial \pbnt}{\partial b}-b\pbnt=0.
\ee
Recall that in the limit $b\to 0$, $(L_+)_0=L_+$ is an elliptic invertible operator in the radial sector and $(L_-)_0=L_-$ is definite positive on $(\mbox{Span}(Q))^{\perp}$ with $Ker(L_-)=\mbox{Span}(Q)$. The following lemma is a standard consequence of the Lax-Milgram theorem and the perturbative theory of uniformly elliptic Schr\"odinger operators, and its proof is left to the reader:

\begin{Lemma}[Invertibility of $(L_+)_b$, $(L_-)_b$]
\label{lemmainversion}
Given $\eta>0$ small, there exists $b^*(\eta)>0$ such that for all $|b|<b^*(\eta)$, the operator $(L_+)_b$ is invertible in the radial sector with 
\be
\label{invertiblitylplusb}
\|e^{(1-\eta)\frac{\theta(|b|r)}{|b|}}(L_+)_b^{-1}f\|_{H^2} \lesssim C_{\eta}\|e^{(1-C\eta)\frac{\theta(|b|r)}{|b|}}f\|_{L^2}.
\ee
Moreover, $(L_-)_b$ admits a lowest eigenvalue $\lambda_b$ with eigenvector $\xi_b$ which are $\mathcal C^1$ functions of $b$ with: 
\be
\label{defdus}
|\lambda_b|+\left\|e^{(1-\eta)\frac{\theta(|b|r)}{|b|}}(\xi_b-Q_p)\right\|_{H^1} \to 0\ \ \mbox{as} \ \ b\to 0
\ee and 
\be
\label{vheioheoghe}
\left|\frac{\partial\lambda_b}{\partial b}\right|+\left\|e^{(1-\eta)\frac{\theta(|b|r)}{|b|}}\frac{\partial \xi_b}{\partial b}\right\|_{L^{\infty}}\leq C \ \ \mbox{as} \ \ b\to 0.
\ee Moreover, 
\be
\label{invertibilitylmoinsb}
\forall f\in (\mbox{Span}(\xi_b))^{\perp}, \ \ \|e^{(1-\eta)\frac{\theta(|b|r)}{|b|}}(L_-)_b^{-1}f\|_{H^2} \lesssim C_{\eta}\|e^{(1-C\eta)\frac{\theta(|b|r)}{|b|}}f\|_{L^{2}}.
\ee
\end{Lemma}

From \fref{invertibilitylmoinsb}, the solvability of the second equation in \fref{eqlpluslmoins} imposes the choice of $\mu_b$:
$$ \left(\mu_b\frac{\partial \pbnt}{\partial b}-b\pbnt,\xi_b\right)=0 \ \ \mbox{i.e.} \ \ \mu_b=\frac{(b\pbnt,\xi_b)}{\left(\frac{\partial \pbnt}{\partial_b},\xi_b\right)}. $$
We now observe the crucial non degeneracy of $\mu_b$ as $b\to 0$ from \fref{convunifqb}, \fref{defdus} and \fref{defrho}: 
\be
\label{defmub}
\mu_b=\frac{(b\pbnt,\xi_b)}{\left(\frac{\partial \pbnt}{\partial_b},\xi_b\right)}\to \frac{|Q_p|_{L^2}^2}{2(\rho,Q_p)}=\frac{4|Q_p|_{L^2}^2}{(1+\sigma_c)|yQ_p|_{L^2}^2} \ \ \mbox{as} \ \ b\to 0,
\ee
where we used the computation from \fref{defrho} and $L_+(\Lambda Q_p)=-2Q_p$:
$$2(\rho,Q_p)=-(L_+\rho,\Lambda Q_p)=-\left(\frac{|y|^2}{4}Q_p,\frac{2}{p-1}Q_p+y\cdot\nabla Q_p\right)=\frac{1+\sigma_c}{4}|yQ_p|_{L^2}^2.$$ The non degeneracy of the denominator as $b\to 0$ and the uniform differentiability properties \fref{convunifqbbis}, \fref{vheioheoghe} ensure that $\mu_b$ is a $\mathcal C^1$ function of $b$ with \be
\label{deifmub}
\left|\frac{\partial \mu_b}{\partial b}\right|\leq C \ \ \mbox{as} \ \ b\to 0.
\ee
Hence from \fref{invertiblitylplusb}, \fref{invertibilitylmoinsb} and the uniform bounds \fref{convunifqb}, \fref{defmub}, \fref{deifmub} we may find $T_b$ solution to \fref{eqlpluslmoins} which is a $\mathcal C^1$ function of $b$ satisfying the uniform bounds \fref{uniformestimatest}.\\

{\bf step 2} Estimate on the error.\\

We now turn to the proof of the estimate of the error \fref{etglobalepsitwo} which amounts estimating the remaining terms in the RHS of \fref{caclulpsitilde}. The nonlinear term is estimated thanks to the homogeneity estimate:
\bee
& & \left|(\pbnt+\sigma_c T_b)|\pbnt+\sigma_c T_b|^{p-1}-(\pbnt)^{p}-\sigma_c p(\pbnt)^{p-1}Re(T_b)-i\sigma_c (\pbnt)^{p-1}Im(T_b)\right|\\
& \lesssim & \left\{\begin{array}{ll} \sigma_c^p|T_p|^p \ \ \mbox{for} \ \ p<2,\\
				\sigma_c^p|T_p|^p+\sigma_c^2|\pbnt|^{p-2}|T_p|^2 \ \ \mbox{for} \ \ p\geq 2,			\end{array} \right .
\eee
and \fref{uniformestimatest} now yields \fref{etglobalepsitwo}, \fref{etglobalepsitwobis}.\\

{\bf step 3} Computation of the invariants.\\

\fref{degenmoment} still holds because $Q_b$ is radially symmetric. \fref{degenviriel}, \fref{qbtsupcritbis}, \fref{propmb} follow from \fref{cnenveo}, \fref{qbtsupcrit}, the decomposition $Q_b=\qbnt+\sigma_c T_be^{-ib\frac{|y|^2}{4}}$, and the fact that Im$(T_b)=O(b)$ by \eqref{eqlpluslmoins}. To compute the energy, we use the Pohozaev multiplier $\Lambda Q_b$ on \fref{defpsib}.  We first integrate by parts to get the general formula: $$Re(\Delta Q_b-Q_b+ib\Lambda Q_b+Q_b|Q_b|^{p-1}, \overline{\Lambda Q_b})=-2E(Q_b)+\sigma_c\left(2E(Q_b)+\int|Q_b|^2\right),$$ and hence from \fref{defpsib}:
\bea
\label{formulaenergyqb}
\nonumber & & 2E(Q_b)-Re\left(\Lambda Q_b,\overline{\Psi_b+i\sigma_c\mu_b\frac{\partial Q_b}{\partial b}}\right)\\
\nonumber & = & 2E(Q_b)-(Re(\Psi_b),\Lambda \S)-(Im(\Psi_b),\Lambda\T)+\sigma_c\mu_b\left[\left(\dtdb,\Lambda\S\right)-\left(\dsdb,\Lambda \T\right)\right]\\
& = & \sigma_c\left(2E(Q_b)+\int|Q_b|^2\right),
\eea
where $\Sigma, \Theta$ are defined by $Q_b=\Sigma+i\Theta$. This together with \fref{qbtsupcritbis} and the estimates on $\Psi_b$ \eqref{preciseest} \fref{etglobalepsitwo} \fref{etglobalepsitwobis} yields the degeneracy \fref{calculhamiltonian}.\\
This concludes the proof of Proposition \ref{propapproximatesolution}.


\subsection{Setting of the bootstrap}


We are now on position to describe the set of initial data $\mathcal O$ leading to the self similar blow up and to set up the bootstrap argument. Pick a Sobolev exponent 
\be
\label{cheioheoh}
\sigma\in\left(\sigma_c,\min\left\{\frac{1}{2},\frac{N}{N+2}\right\}\right)
\ee
 independent of $p$ and close enough to $0$.

\begin{Definition}[Geometrical description of the set $\mathcal O$]
\label{defdefp}
Pick a number $\nu_0>0$ small enough. Then for $p\in(p_c,p^*(\nu_0))$ with $p^*(\nu_0)$ close enough to $p_c$, we let $\mathcal O$ be the set of initial data $u_0\in H^1$ of form:
 $$u_0(x)=\frac{1}{\lambda_0^{\frac{2}{p-1}}}(Q_{b_0}+\e_0)\left(\frac{x-x_0}{\lambda_0}\right)e^{i\gamma_0}$$ for some $(\lambda_0,b_0,x_0,\gamma_0)\in \RR^*_+\times \RR^*_+\times \RR^N\times \RR$ with the following controls:\\
(i) $b_0$ is in the self similar asymptotics \fref{sigmac}:
\be
\label{comparaisonbalphainit}
\Gamma_{b_0}^{1+\nu^{10}_0}\leq \sigma_c \leq \Gamma_{b_0}^{1-\nu^{10}_0};
\ee
(ii) Smallness of the scaling parameter: 
\be
\label{lambdasmallimit}
0<\lambda_0\leq \Gamma_{b_0}^{100};
\ee
(iii) Degeneracy of the energy and the momentum:
\be
\label{initienergyomentum}
\lambda_0^{2(1-
\sigma_c)}|E_0|+\lambda_0^{1-2\sigma_c}\left|Im\left(\int\nabla u_0\overline{u}_0\right)\right|<\Gamma_{b_0}^{50};
\ee
(iv) $\dot{H}^1\cap \dot{H}^{\sigma}$ smallness of the excess of $L^2$ mass:
\be
\label{esteinit}
\int||\nabla|^{\sigma}\e_0|^2+\int|\nabla \e_0|^2+\int |\e_0|^2e^{-|y|}< \Gamma_{b_0}^{1-\nu_0}.
\ee
\end{Definition}

\begin{remark} Hidden in Definition \ref{defdefp} is the choice of the parameter $\eta>0$ entering in the 
construction of $Q_b^{(0)}$ in Proposition \ref{profilespbzero}. In all what follows, we will need the fact that given a universal constant $C>0$, we have the control $$\gab^{1-C\eta}\leq \gab^{1-\nu_0^{50}}.$$ This holds provided $\eta>0$ has been chosen small enough with respect to $\nu_0$. We may for example take $$\eta=\nu_0^{100}.$$ Now, pushing $\eta\to 0$ requires pushing $b\to 0$ or equivalently $p\to p_c$ from \fref{comparaisonbalphainit}. This is how the asymptotics \fref{loisigma} follows.
\end{remark}

\begin{remark} Observe that the set $\mathcal O$ is a non empty open set in $H^1$. Indeed, pick a small parameter $\nu_0$ and $p$ close enough to $p_c$. Pick $b_0>0$ such that \fref{comparaisonbalphainit} holds, and pick then $\lambda_0>0$ such that \fref{lambdasmallimit} holds. Let $f(y)$ be smooth, real, radial and compactly supported in the ball $|y|\leq 1$ and such that $(f,Q)=1$, and let $\e_0=\mu_0f$ with $\mu_0$ to be chosen. From $$\frac{d}{d\mu_0}E(Q_{b_0}+\mu_0f)_{|\mu_0=0}=-(f,Q)(1+o(1))=-1+o(1) \ \ \mbox{as} \ \ b_0\to 0,$$ and the degeneracy of the Hamiltonian \fref{calculhamiltonian}, we may find $\mu_0=O(\Gamma_{b_0}^{1-2\nu_0})$ such that $|E(Q_{b_0}+\e_0)|\lesssim \Gamma_{b_0}^{100}$ so that \fref{initienergyomentum} holds. \fref{esteinit} now follows from the size of $\mu_0$. 
\end{remark}

Let now $u_0\in \mathcal O$ and $u(t)$ be the corresponding solution to \fref{nls} with maximum life time interval $[0,T)$, $0<T\leq +\infty$. Using the regularity $u\in \mathcal C([0,T),H^1)$ and standard modulation theory, we can find a small interval $[0,T^*)$ such that for all $t\in [0,T^*)$, $u(t)$ admits a unique geometrical decomposition \
\be
\label{decompou}
u(t,x)=\frac{1}{\lambda^{\frac{2}{p-1}}(t)}(Q_{b(t)}+\e)\left(t,\frac{x-x(t)}{\lambda(t)}\right)e^{i\gamma(t)}
\ee
 where uniqueness follows from the freezing of orthogonality conditions:  $\forall t\in[0,T^*]$,
\be
\label{orthe1b}
\left(\e_1(t),|y|^2\S\right)+\left(\e_2(t),|y|^2\T\right)=0 ,
\ee
\be
\label{orthe2b}
\left(\e_1(t),y\S\right)+\left(\e_2(t),y\T\right)=0 ,
\ee
\be
\label{orthe3b}
\left(\e_2(t),\Lambda^2\S\right)-\left(\e_1(t),\Lambda^2\T\right)=0 ,
\ee
\be
\label{orthe4b}
\left(\e_2(t),\Lambda \S\right)-\left(\e_1(t),\Lambda \T\right)=0,
\ee
where we have denoted: 
$$
\e=\e_1+i\e_2, \ \ Q_b=\Sigma+i\Theta$$ in terms of real and imaginary parts. See \cite{MR1}, \cite{MR2} for related statements. Moreover, the parameters $(\lambda(t),b(t),x(t),\gamma(t))\in \RR^*_+\times \RR^*_+\times\RR^N\times\RR$ are $\mathcal C^1$ functions of time and $\e\in \mathcal C([0,T^*),H^1)$ with a priori bounds: $\forall t\in [0,T^*)$,
\be
\label{comparaisonbalphaboot}
\Gamma_{b(t)}^{1+\nu^2_0}\leq \sigma_c \leq \Gamma_{b(t)}^{1-\nu^2_0},
\ee
\be
\label{estlambdaboot}
0<\lambda(t)\leq \Gamma_{b(t)}^{10},
\ee
\be
\label{energyomentumboot}
[\lambda(t)]^{2(1-\sigma_c)}|E_0|+[\lambda(t)]^{1-2\sigma_c}\left|Im\left(\int\nabla u_0\overline{u}_0\right)\right|<\Gamma_{b(t)}^{10},
\ee
\be
\label{estenboot}
\int|\nabla \e(t)|^2+\int |\e(t)|^2e^{-|y|}<\Gamma_{b(t)}^{1-20\nu_0},
\ee
\be
\label{esteisigmaboot}
\int||\nabla|^{\sigma}\e(t)|^2\leq \Gamma_{b(t)}^{1-50\nu_0}.
\ee

\begin{remark} The strict $\dot{H}^1$ subcriticality of the problem implies: $$\sigma_c=\frac{N}{2}-\frac{2}{p-1}<\frac{N}{2}-\frac{N}{p+1}.$$ Hence from Sobolev embedding, the $\dot{H}^{\sigma}$ control \fref{esteisigmaboot} together with \eqref{estenboot} ensures for $\sigma$ close enough to $\sigma_c$:
\be
\label{estnonlinearterm}
\int|\e|^{p+1}\lesssim |\e|_{\dot{H}^{\frac{N}{2}-\frac{N}{p+1}}}^{p+1}\leq \left[\Gamma_{b}^{1-50\nu_0}\right]^{\frac{p+1}{2}}\lesssim \gab^{1+z_0}
\ee 
for some universal constant $z_0>0$ independent of $p$ for $p$ close enough to $p_c$.
\end{remark}

Our main claim is that the above regime is a trapped regime:

\begin{Prop}[Bootstrap]
\label{propboot}
There holds: $\forall t\in [0,T^*)$,
\be
\label{comparaisonbalphabootalpha}
\Gamma_{b(t)}^{1+\nu^4_0}\leq \sigma_c \leq \Gamma_{b(t)}^{1-\nu^4_0},
\ee
\be
\label{estlambdabootfinal}
0<\lambda(t)\leq \Gamma_{b(t)}^{20},
\ee
\be
\label{energyomentumfinal}
[\lambda(t)]^{2(1-\sigma_c)}|E_0|+[\lambda(t)]^{1-2\sigma_c}\left|Im\left(\int\nabla u_0\overline{u}_0\right)\right|<\Gamma_{b(t)}^{20},
\ee
\be
\label{estebootfinal}
\int|\nabla \e(t)|^2+\int |\e(t)|^2e^{-|y|}<\Gamma_{b(t)}^{1-10\nu_0},
\ee
\be
\label{esteisigmafinal}
\int||\nabla|^{\sigma}\e(t)|^2\leq \Gamma_{b(t)}^{1-45\nu_0}
\ee
and hence $$T^*=T.$$
\end{Prop}

The next section is devoted to the derivation of the key dynamical controls at the heart of the proof of the bootstrap Proposition \ref{propboot} which is proved in section \ref{sectionproofthm}. Theorem \ref{thmmain} will be a simple consequence of Proposition \ref{propboot}.


\section{Control of the self similar dynamics}
\label{sectionthree}


In this section, we exhibit the two Lyapounov type functionals which will allow us to lock the self similar dynamics and prove the bootstrap Proposition \ref{propboot}. The key is the dynamical lock \fref{comparaisonbalphabootalpha} of the geometrical parameter $b(t)$ which controls the selfsimilarity of blow up from the modulation equation $$b\sim-\lsl=-\lambda\lambda_t.$$ This corresponds to an upper bound and a lower bound of $b$. The proof will follow by somehow bifurcating from the log-log analysis and by tracking the effect of the leading order $\sigma_c$ deformation in the $Q_b$ profile. The lower bound on b is a consequence of the local virial estimate type of control first derived in \cite{MR2}, see Proposition \ref{propviriel}, the upper bound follows from the sharp log-log analysis derived in \cite{MR4} and uses very strongly the $L^2$ conservation law, see Proposition \ref{lemmadiffineqb}. Remember also that we have no a priori orbital stability bound on neither $b$ or $\e$, and indeed the upper {\it pointwise} control \fref{estebootfinal} on $\e$ requires {\it both} monotonicity properties, see step 1 of the proof of Proposition \ref{propboot} in section \ref{sectionclosingboot}.


\subsection{Preliminary estimates on the decomposition}


Let us recall the geometrical decomposition \fref{decompou}:
$$
u(t,x)=\frac{1}{\lambda^{\frac{2}{p-1}}(t)}(Q_{b(t)}+\e)\left(t,\frac{x-x(t)}{\lambda(t)}\right)e^{i\gamma(t)}
$$
and derive the modulation equations on the geometrical parameters and preliminary estimates inherited from the conservation laws. We let the rescaled time $$s(t)=\int_0^t\frac{d\tau}{\lambda^2(\tau)},\ \ s^*=s(T^*)\in(0,+\infty],$$ and compute the equation of $\e$ in terms of real and imaginary parts on $[0,s^*)$:
\bea
\label{partiereelleb}
\nonumber(b_s-\sigma_c\mu_b)\dsdb+\partial_s\er-M_{-}(\e)+b\Lambda\e_1 & =& \left(\lsl+b\right)\Lambda \S+\tgamma_s\T+\xsl\cdot\nabla\S\\
\nonumber & + & \left(\lsl+b\right)\Lambda\e_1+\tgamma_s\e_2+\xsl\cdot\nabla\e_1\\
& + & Im(\Psi_b)-R_2(\e)
\eea
\bea
\label{partieimaginaireb}
\nonumber(b_s-\sigma_c\mu_b)\dtdb+\partial_s\e_2+M_+(\e)+b\Lambda\e_2& = & \left(\lsl+b\right)\Lambda \T-\tgamma_s\S+\xsl\cdot\nabla\T\\
\nonumber& + & \left(\lsl+b\right)\Lambda\e_2-\tgamma_s\e_1+\xsl\cdot\nabla\e_2\\
 & - & Re(\Psi_b)+R_1(\e),
\eea
with $\tilde{\gamma}(s)=-s+\gamma(s)$ and $\Psi_b$ given by \fref{defpsib}. Here we also denoted $M=(M_+,M_-)$  the linear operator close to $\qb$, explicitly: $$M_+(\e)=-\Delta \e_1+\e_1-\left(1+(p-1)\frac{\S^2}{|\qb|^2}\right)|\qb|^{p-1}\e_1-(p-1)\S\T|\qb|^{p-3}\e_2,$$$$M_-(\e)=-\Delta \e_2+\e_2-\left(1+(p-1)\frac{\T^2}{|\qb|^2}\right)|\qb|^{p-1}\e_2-(p-1)\S\T|\qb|^{p-3}\e_1.$$ The non linear interaction terms are explicitly: 
\bea
\label{run}
R_1(\e) & = & (\e_1+\S)|\e+\qb|^{p-1}-\S|\qb|^{p-1}\\
\nonumber & - & \left(1+(p-1)\frac{\S^2}{|\qb|^2}\right)|\qb|^{p-1}\e_1-(p-1)\S\T|\qb|^{p-3}\e_2,
\eea
\bea
\label{rdeux}
R_2(\e) & = & (\e_2+\T)|\e+\qb|^{p-1}-\T|\qb|^{p-1}\\
\nonumber & - & \left(1+(p-1)\frac{\T^2}{|\qb|^2}\right)|\qb|^{p-1}\e_2-(p-1)\S\T|\qb|^{p-3}\e_1.
\eea
We now claim the following preliminary estimates on the decomposition:

\begin{Lemma}
\label{propdecompbun}
 There holds for some universal constants $C>0$, $\delta(p)\to 0$ as $p\to p_c$ and for all $s\in [0,s^*)$:\\
(i)Estimates induced by the conservation of energy and momentum:
\be
\label{controlelenergyb}
 |2(\e_1,\S)+2(\e_2,\T)|\leq C\left(\int |\nabla \e|^2+\normeldloc\right)+\gab^{1-11\nu_0}
\ee
\be
\label{degenracymometn}
|(\e_2,\nabla\S)|\leq \delta(p)\left(\int |\nabla\e |^2+\normeldloc\right)^{\frac{1}{2}}+\gab^{1-50\nu_0}.
\ee
(ii) Estimates on the modulation parameters:
\be
\label{smalnesscoeffpairb}
\left|\lsl+b\right|+|b_s|\leq C\left(\int|\nabla \e|^2+\int |\e| ^2e^{-|y|}\right)+\gab^{1-11\nu_0},
\ee
\bea
\label{smallnesscoeffimpairb}
\nonumber \left|\tgamma_s-\frac{1}{|DQ|_2^2}(\e_1,L_+D^2Q)\right|+\left|\xsl\right| & \leq & \delta(p)\left(\int|\nabla\e|^2e^{-2(1-\eta)\frac{\theta(b|y|)}{b}}+\int |\e| ^2e^{-|y|}\right)^{\frac{1}{2}}\\
& + & C\int|\nabla\e|^2+\gab^{1-11\nu_0}.
\eea
\end{Lemma}

{\bf Proof of Lemma \ref{propdecompbun}}\\

It relies as in \cite{MR4} on the expansion of the conservation laws, the choice of orthogonality conditions for $\e$ and the bootstrapped controls \fref{comparaisonbalphaboot}, \fref{estlambdaboot}, \fref{energyomentumboot}, \fref{estenboot}, \fref{esteisigmaboot}.\\

{\bf step1 } Expansion of the conservation laws.\\

\fref{controlelenergyb} and \fref{degenracymometn} follow from the expansion of the momentum and the energy using the decomposition \fref{decompou} and the estimates of Proposition \ref{propapproximatesolution}.\\
For the momentum: 
$$
2Im(\e,\overline{\nabla Q_b})=Im\int\nabla \e\overline{\e}-\lambda^{1-2\sigma_c}Im\left(\int\nabla u_0\overline{u_0}\right)
$$
from which using \fref{energyomentumboot}:
\bee
 |(\e_2,\nabla Q)| & \leq & 
 \delta(p)\left(\int |\nabla\e |^2+\normeldloc\right)^{\frac{1}{2}}+C|\e|_{\dot{H}^{\frac{1}{2}}}^2+\gab^{10}\\
 & \leq & \delta(p)\left(\int |\nabla\e |^2+\normeldloc\right)^{\frac{1}{2}}+\gab^{1-50\nu_0}
 \eee
 where we interpolated $\dot{H}^{\frac{1}{2}}$ between $\dot{H}^{\sigma}$ and $\dot{H}^1$ and used \fref{estenboot}, \fref{esteisigmaboot}.\\
For the energy:
\bea
\label{energyconservation}
\nonumber & & 2\left(\e_1,\Sigma+b\Lambda\T-Re(\Psi_b)+\sigma_c\mu_b\dtdb\right)+2\left(\e_2,\T-b\Lambda \Sigma-Im(\Psi_b)-\sigma_c\mu_b\dsdb\right)\\
\nonumber & = & -2\lambda^{2(1-\sigma_c)}E_0+2E(\qb)+ (M_+(\e),\e_1)+(M_-(\e),\e_2)-\int|\e|^2\\
& - & \frac{2}{p+1}\int F(\e)
\eea
where $F(\e)$ is the formally cubic part of the potential energy: 
\bea
\label{deffe}
\nonumber F(\e)& = & |Q_b+\e|^{p+1}-|Q_b|^{p+1}-(p+1)Re(\e,\overline{Q_b|Q_b|^{p-1}})\\
\nonumber & - &  (p+1)\left(1+(p-1)\frac{\S^2}{|\qb|^2}\right)|\qb|^{p-1}\e^2_1-(p+1)\left(1+(p-1)\frac{\T^2}{|\qb|^2}\right)|\qb|^{p-1}\e_2^2\\
& - & 2(p+1)(p-1)\S\T|\qb|^{p-3}\e_2\e_1.
\eea
We then use standard homogeneous estimates and Sobolev embeddings like in \cite{MR1}, \cite{MR2} to estimate the nonlinear term $F(\e)$. Indeed, choose $\mu$ such that $0<\mu<\min(1,p-1)$, then by homogeneity: $$|F(\e)|\lesssim |\e|^{p+1}+|Q_b|^{p-1-\mu}|\e|^{2+\mu},$$ and thus from Holder with $\frac{1}{2+\mu}=\frac{\alpha}{2}+\frac{1-\alpha}{p+1}$ and the bootstrap controls \fref{estenboot}, \fref{estnonlinearterm}:
\bea
\label{estfe}
\nonumber 
\int |F(\e)| & \lesssim & \int |\e|^{p+1}+\left|\e e^{-C|y|}\right|_{L^2}^{\alpha(2+\mu)}|\e|_{L^{p+1}}^{(1-\alpha)(2+\mu)}\\
& \lesssim & \gab^{1+z_0}+\left(\gab^{1-20\nu_0}\right)^{2+\mu}\lesssim \gab^{1+z_0}
\eea for some $z_0>0$ independent of $p$.
Injecting this into \fref{energyconservation} together with the degeneracy estimate \fref{calculhamiltonian}, the bootstrap bound \fref{comparaisonbalphaboot}, \fref{estlambdaboot} and the orthogonality condition \fref{orthe4b} yields \fref{controlelenergyb}.

\begin{remark}
\label{rkenergy} Note that the algebraic formula \fref{energyconservation} does not use the orthogonality conditions on $\e$.
\end{remark}

{\bf step 2} Computation of the modulation parameters.\\

The estimates \fref{smalnesscoeffpairb}, \fref{smallnesscoeffimpairb} are obtained by computing the geometrical parameters from the choice of orthogonality conditions \fref{orthe1b}, \fref{orthe2b}, \fref{orthe3b}, \fref{orthe4b} and by relying on the estimates induced by the conservation of energy and momentum \eqref{controlelenergyb} \eqref{degenracymometn}. The computation is the same like the one performed in \cite{MR4} up to $O(\sigma_c)=O(p-p_c)$ terms which are brutally estimated in absolute value using the bootstrap bound \fref{comparaisonbalphaboot}. The detail of this is left to the reader.\\ 
This concludes the proof of Lemma \ref{propdecompbun}.


\subsection{Local virial identity}


We now proceed through the derivation of the local virial control. The corresponding monotonicity property was first discovered in \cite{MR1} and will yield a strict {\it lower bound} on $b$.

\begin{Prop}[Local virial identity]
\label{propviriel}
There holds for some universal constant $c_1>0$ the lower bound:
\be
\label{lowerboundbs}
\forall s\in[0,s^*), \ \ b_s\geq c_1\left(\sigma_c+\int|\nabla \e|^2+\int|\e|^2e^{-|y|}-\gab^{1-\nu^6_0}\right).\ee
\end{Prop}

\begin{remark} The super critical effect relies in the presence of the term $\sigma_c$ in the RHS of \fref{lowerboundbs}. The positive sign is crucial and structural and is the reason why we had to push the construction of $Q_b$ to an order $\sigma_c^2$.
\end{remark}

{\bf Proof of Proposition \ref{propviriel}}\\

{\bf step1} Algebraic derivation of the $b_s$ law.\\

The first step is a careful derivation of the modulation equations for b. For further use, we shall exhibit first a general formula which {\it does not rely on the specific choice of orthogonality} conditions in $\e$. Once the formula is derived, the use of the orthogonality conditions and suitable coercivity properties inherited from the Spectral Property stated in the introduction will yield the claim.\\

Take the inner product of \fref{partiereelleb} with $(-\Lambda\T)$ and \fref{partieimaginaireb} with $\Lambda \S$ and sum the obtained identities using \fref{integraoarts} to get:
\bea
\label{firstexpressionbs}
\nonumber  & & b_s\left(\left(\dtdb,\Lambda\S\right)-\left(\dsdb,\Lambda\T\right)-\left(\e_2,\Lambda\dsdb\right)+\left(\e_1,\Lambda \dtdb\right)\right)+\left\{(\e_2,\Lambda\S)-(\e_1,\Lambda\T)\right\}_s= \\
 & - & (M_+(\e)+b\Lambda \e_2,\Lambda\S)-(M_-(\e)-b\Lambda\e_1,\Lambda\T) - \tgamma_s\left\{(\e_1,\Lambda\S)+(\e_2,\Lambda\T)\right\}\\
 \nonumber & - &  \left(\lsl+b\right)\left\{(\e_2,\Lambda^2\S+2\sigma_c\Lambda\S)-(\e_1,\Lambda^2\T+2\sigma_c\Lambda\T)\right\}-  \xsl\cdot\left\{(\e_2,\nabla\Lambda\S)-(\e_1,\nabla\Lambda \T)\right\}\\
\nonumber & - & Re\left(\Lambda Q_b,\overline{\Psi_b+i\sigma_c\mu_b\frac{\partial Q_b}{\partial b}}\right)+\sigma_c\tgamma_s|Q_b|_{L^2}^2+(R_1(\e),\Lambda \Sigma)+(R_2(\e),\Lambda\T).
\eea
We compute the linear term in $\e$ in \fref{firstexpressionbs}. For this, consider the $Q_b$ equation \fref{defpsib}, compute the equation satisfied by $\mu^{\frac{2}{p-1}}Q(\mu y)$ and differentiate the obtained identity with respect to $\mu$ at $\mu=1$. This yields:
$$M_+(\Lambda Q_b)+b\Lambda^2\T=-2\left[\S+b\Lambda\T-Re(\Psi_b)+\sigma_c\mu_b\dtdb\right]+Re(\Lambda\Psi_b)-\sigma_c\mu_b\Lambda\dtdb,$$ $$M_-(\Lambda Q_b)-b\Lambda^2\Sigma=-2\left[\T-b\Lambda \S-Im(\Psi_b)-\sigma_c\mu_b\dsdb\right]+Im(\Lambda\Psi_b)+\sigma_c\mu_b\Lambda \dsdb.$$ Integrating by parts and using the conservation of the energy \fref{energyconservation}, we obtain the following identity:
\bee
& & - (M_+(\e)+b\Lambda \e_2,\Lambda\S)-(M_-(\e)-b\Lambda\e_1,\Lambda\T)\\
& = & -(\e_1,M_+(\Lambda Q_b)+b\Lambda^2\T+2b\sigma_c\Lambda\T)-(\e_2,M_-(\Lambda Q_b)-b\Lambda^2\S-2b\sigma_c\S)\\
& = & \left(\e_1, 2\left[\S+b\Lambda\T-Re(\Psi_b)+\sigma_c\mu_b\dtdb\right]-Re(\Lambda\Psi_b)+\sigma_c\mu_b\Lambda\dtdb-2b\sigma_c\Lambda\T\right)\\
& + & \left(\e_2,2\left[\T-b\Lambda \S-Im(\Psi_b)-\sigma_c\mu_b\dsdb\right]-Im(\Lambda\Psi_b)-\sigma_c\mu_b\Lambda \dsdb+2b\sigma_c\Lambda\S\right)\\
& = & 2b\sigma_c\left[-(\e_1,\Lambda\T)+(\e_2,\Lambda\S)\right]-(\e_1,Re(\Lambda\Psi_b))-(\e_2,Im(\Lambda\Psi_b))\\
& + & \sigma_c\mu_b\left[\left(\e_1,\Lambda\dtdb\right)-\left(\e_2,\Lambda\dsdb\right)\right]- 2\lambda^{2(1-\sigma_c)}E_0+2E(\qb)\\
& + & (M_+(\e),\e_1)+(M_-(\e),\e_2)-\int|\e|^2-\frac{2}{p+1}\int F(\e)\\
& = & \delta(p)\sigma_c-(\e_1,Re(\Lambda\Psi_b))-(\e_2,Im(\Lambda\Psi_b))- 2\lambda^{2(1-\sigma_c)}E_0+2E(\qb)\\
& + & (M_+(\e),\e_1)+(M_-(\e),\e_2)-\int|\e|^2-\frac{2}{p+1}\int F(\e)
\eee
where we recall that $\delta(p)$ denotes a generic constant $\delta(p)\to 0$ as $p\to p_c$. We now inject this into \fref{firstexpressionbs}. A key here is to use the Pohozaev identity \fref{formulaenergyqb}:
$$2E(Q_b)-Re\left(\Lambda Q_b,\overline{\Psi_b+i\sigma_c\mu_b\frac{\partial Q_b}{\partial b}}\right)= \sigma_c\left(2E(Q_b)+\int|Q_b|^2\right)
$$ to generate a {\it nonnegative term} in the RHS of \fref{firstexpressionbs}:
\bea
\label{firstexpressionbsbisbisbis}
\nonumber  & & b_s\left(\left(\dtdb,\Lambda\S\right)-\left(\dsdb,\Lambda\T\right)-\left(\e_2,\Lambda\dsdb\right)+\left(\e_1,\Lambda \dtdb\right)\right)+\left\{(\e_2,\Lambda\S)-(\e_1,\Lambda\T)\right\}_s=\\
\nonumber  & & \delta(p)\sigma_c+\sigma_c(2E(Q_b)+|Q_b|_{L^2}^2+\tgamma_s|Q_b|_{L^2}^2)-2\lambda^{2(1-\sigma_c)}E_0\\
& - & (\e_1,Re(\Lambda\Psi_b))-(\e_2,Im(\Lambda\Psi_b))\\
& + &  (M_+(\e),\e_1)+(M_-(\e),\e_2)-\int|\e|^2+(R_1(\e),\Lambda \Sigma)+(R_2(\e),\Lambda\T)\\
 \nonumber & - & \tgamma_s\left\{(\e_1,\Lambda\S)+(\e_2,\Lambda\T)\right\}- \xsl\cdot\left\{(\e_2,\nabla\Lambda\S)-(\e_1,\nabla\Lambda \T)\right\}\\
\nonumber  & - &  \left(\lsl+b\right)\left\{(\e_2,\Lambda^2\S+2\sigma_c\Lambda\S)-(\e_1,\Lambda^2\T+2\sigma_c\Lambda\T)\right\}-\frac{2}{p+1}\int F(\e).
\eea
It remains to extract the formally cubic term in $\e$ in the RHS \fref{firstexpressionbsbisbisbis}. We let:
\bea
\label{gone}
\nonumber G_1(\e) & = & R_1(\e)-\frac{p-1}{2}|Q_b|^{p-5}(p\S^3+3\S\T^2)\e_1^2-\frac{p-1}{2}|Q_b|^{p-5}(\S^3+(p-2)\S\T^2)\e_2^2\\
& - & (p-1)|Q_b|^{p-5}(\T^3+(p-2)\S^2\T)\e_1\e_2,
\eea
\bea
\label{godeux}
\nonumber G_2(\e) & = & R_2(\e)-\frac{p-1}{2}|Q_b|^{p-5}(\T^3+(p-2)\T\S^2)\e_1^2-\frac{p-1}{2}|Q_b|^{p-5}(p\T^3+3\T\S^2)\e_2^2\\
& - & (p-1)|Q_b|^{p-5}(\S^3+(p-2)\T^2\S)\e_1\e_2,
\eea
and eventually arrive at the following algebraic virial identity:
\bea
\label{firstexpressionbsbisbis}
\nonumber  & & b_s\left(\left(\dtdb,\Lambda\S\right)-\left(\dsdb,\Lambda\T\right)-\left(\e_2,\Lambda\dsdb\right)+\left(\e_1,\Lambda \dtdb\right)\right)+\left\{(\e_2,\Lambda\S)-(\e_1,\Lambda\T)\right\}_s=\\
\nonumber  & & \delta(p)\sigma_c+\sigma_c(2E(Q_b)+|Q_b|_{L^2}^2+\tgamma_s|Q_b|_{L^2}^2)-2\lambda^{2(1-\sigma_c)}E_0\\
& - & (\e_1,Re(\Lambda\Psi_b))-(\e_2,Im(\Lambda\Psi_b))+H_p(\e,\e)-\tgamma_s\left\{(\e_1,\Lambda\S)+(\e_2,\Lambda\T)\right\}\\
 \nonumber & - & \xsl\cdot\left\{(\e_2,\nabla\Lambda\S)-(\e_1,\nabla\Lambda \T)\right\}+ E(\e,\e)+(G_1(\e),\Lambda \Sigma)+(G_2(\e),\Lambda\T)\\
\nonumber  & - & \left(\lsl+b\right)\left\{(\e_2,\Lambda^2\S+2\sigma_c\Lambda\S)-(\e_1,\Lambda^2\T+2\sigma_c\Lambda\T)\right\}-\frac{2}{p+1}\int F(\e),
\eea
where the virial quadratic form $H_p$ can be expressed in the form:
\bea
\label{formequadravirielp}
H_p(\e,\e) & = & \int|\nabla\e|^2+\frac{p(p-1)}{2}\int y\cdot\nabla Q_pQ_p^{p-2}\e_1^2\\
\nonumber & + & \frac{(p-1)}{2}\int y\cdot\nabla Q_pQ_p^{p-2}\e_2^2
\eea
and the error $E(\e,\e)$ simply comes from the error made by replacing $Q_b$ by $Q$ in the potential terms and is easily estimated thanks to Proposition \ref{propapproximatesolution} by:
\be
\label{esterrore}
|E(\e,\e)|\leq \delta(p)\left(\int|\nabla \e|^2+\normeldloc\right)
\ee
with $\delta(p)\to 0$ as $p\to p_c.$ \\

{\bf step 2} Estimates of nonlinear terms and coercivity of the quadratic form.\\

Observe from the construction of $Q_b$ that 
\bea
\label{calculprelim}
\nonumber& &  \left(\dtdb,\Lambda\S\right)-\left(\dsdb,\Lambda\T\right)  =  -Im\left(\Lambda Q_b,\overline{\frac{\partial Q_b}{\partial b}}\right)\\
& = & -Im\left(\Lambda P_b-ib\frac{|y|^2}{2}P_b,\overline{\frac{\partial P_b}{\partial b}-i\frac{|y|^2}{4}P_b}\right)
=\frac{|yQ_{p_c}|_{L^2}^2}{4}+O(|b|+\sigma_c).
\eea
We then inject the orthogonality conditions \fref{orthe1b}, \fref{orthe2b}, \fref{orthe3b}, \fref{orthe4b} and \fref{calculprelim} into \fref{firstexpressionbsbisbis} to derive the algebraic identity:
\bea
\label{firstexpressionbsbis}
 & & \frac{|yQ_{p_c}|_{L^2}^2}{4}(1+\delta(p))b_s\\
\nonumber  & = & \delta(p)\sigma_c+\sigma_c(2E(Q_b)+|Q_b|_{L^2}^2+\tgamma_s|Q_b|_{L^2}^2)-2\lambda^{2(1-\sigma_c)}E_0-(\e_1,Re(\Lambda\Psi_b))-(\e_2,Im(\Lambda\Psi_b))\\
\nonumber & + &  H_p(\e,\e)-\tgamma_s\left\{(\e_1,\Lambda\S)+(\e_2,\Lambda\T)\right\}-  \xsl\cdot\left\{(\e_2,\nabla\Lambda\S)-(\e_1,\nabla\Lambda \T)\right\}\\
\nonumber & +&  E(\e,\e)+(G_1(\e),\Lambda \Sigma)+(G_2(\e),\Lambda\T)-\frac{2}{p+1}\int F(\e)
\eea
with $\delta(p)\to 0$ as $p\to p_c$.
Let us now estimate all the terms in the RHS of \fref{firstexpressionbsbis}. First observe from \fref{qbtsupcritbis}, \fref{calculhamiltonian} and the smallness of $\e$ the non degeneracy:
\be
\label{estalphalower}
\delta(p)\sigma_c+\sigma_c(2E(Q_b)+|Q_b|_{L^2}^2+\tgamma_s|Q_b|_{L^2}^2)\geq \frac{\sigma_c|Q_{p_c}|_{L^2}^2}{2}.
\ee
The nonlinear terms are estimated as in \cite{MR2}, \cite{MR4} using homogeneity estimates, Sobolev embeddings and the bootstrap bounds \fref{estenboot}, \fref{estnonlinearterm} as for the proof of \fref{estfe}:
\bea
\label{estnolinea}
\nonumber & & \left|E(\e,\e)+(G_1(\e),\Lambda \Sigma)+(G_2(\e),\Lambda\T)-\frac{2}{p+1}\int F(\e)\right|\\
\nonumber & \leq & \delta(p)\left(\int|\nabla \e|^2+\normeldloc\right)+\gab^{1+z_0}.
\eea
We now focus on the quadratic terms in \fref{firstexpressionbsbis}. We inject the estimates of the geometrical parameters \fref{smalnesscoeffpairb}, \fref{smallnesscoeffimpairb}. Using the simple bound  $$\left|e^{c|y|}\left[|Q_{p}-Q_{p_c}|+|\nabla Q_p-\nabla Q_{p_c}|\right]\right|_{L^{\infty}}\to 0 \ \ \mbox{as} \ \ p\to p_c,$$ we may view the obtained quadratic form together with $H_p$ as a perturbation of the one obtained in the $L^2$ critical case:
 \bee
 & &  H_p(\e,\e)-\tgamma_s\left\{(\e_1,\Lambda\S)+(\e_2,\Lambda\T)\right\}- \xsl\cdot\left\{(\e_2,\nabla\Lambda\S)-(\e_1,\nabla\Lambda \T)\right\}\\
 &= & \tilde{H}(\e,\e)+E_1(\e,\e)
 \eee
with 
\bee
\tilde{H}(\e,\e) & = & \int |\nabla \e|^2+\frac {2}{N}\left(1+\frac{4}{N}\right)\int Q_{p_c}^{\frac{4}{N}-1}y\cdot\nabla Q_{p_c} \e_1^2+\frac {2}{N}\int Q_{p_c}^{\frac{4}{N}-1}y\cdot\nabla Q_{p_c} \e_1^2\\
& - &\frac{1}{|DQ_{p_c}|_{L^2}^2}(\e_1,L_+D^2Q_{p_c})(\e_1,DQ_{p_c}),
\eee
$$|E_1(\e,\e)|\leq \delta(p)\left(\int|\nabla \e|^2+\normeldloc\right)+\gab^{2-30\nu_0}.$$
We now recall from \cite{MR3} the following coercivity property which is a consequence of the spectral property stated in the introduction: $\forall \e=\e_1+i\e_2\in H^1$, 
\bee
\tilde{H}(\e,\e)& \geq&   c_0\left(\int|\nabla \e|^2+\normeldloc\right)\\
\nonumber & - & \left\{(\e_1,Q_{p_c})^2+(\e_1,|y|^2Q_{p_c})^2+(\e_1,yQ_{p_c})^2+(\e_2,DQ_{p_c})^2+(\e_2,D^2Q_{p_c})^2+(\e_2,\nabla Q_{p_c})^2\right\}
\eee
for some universal constant $c_0>0$, and hence our choice of orthogonality conditions, \fref{genrtor} and the degeneracy estimates \fref{controlelenergyb}, \eqref{degenracymometn} ensure:
\be
\label{vneonveon}
\tilde{H}(\e,\e)\geq  \frac{c_0}{2}\left(\int|\nabla \e|^2+\normeldloc\right)-\gab^{\frac{3}{2}}
\ee
for $\nu_0=\nu_0(p)>0$ small enough. We now inject \fref{calculprelim}, the orthogonality condition \fref{orthe4b}, the nondegeneracy estimate \fref{estalphalower}, the estimates on nonlinear terms \fref{estnolinea} and the coercivity property \fref{vneonveon} into \fref{firstexpressionbsbis} to derive
\bee
b_s & \geq & \frac{c_0}{4}\left(\sigma_c +\int|\nabla \e|^2+\normeldloc\right)-\gab^{1-\nu_0^7}\\
& - & |(\e_1,Re(\Lambda\Psi_b))+(\e_2,Im(\Lambda\Psi_b))|.
\eee
A crude bound for the remaining linear term is derived from \eqref{preciseest} \fref{etglobalepsitwo} \fref{etglobalepsitwobis}:
\be
\label{remaingicrudebound}
|(\e_1,Re(\Lambda\Psi_b)+(\e_2,Im(\Lambda\Psi_b)|\leq \gab^{1-C\eta}+\frac{c_0}{10}\left(\sigma_c +\int|\nabla \e|^2+\normeldloc\right),
\ee and \fref{lowerboundbs} follows for $\eta<\eta(\nu_0)$ chosen small enough.\\
This concludes the proof of Proposition \ref{propviriel}.


\subsection{Refined local virial identity and introduction of the radiation}

We now proceed through a refinement of the local virial estimate \fref{lowerboundbs} and adapt the analysis in \cite{MR4}.\\
We start with introducing a localized version of the radiation $\zeta_b$ introduced in Lemma \ref{outgoingradiation} due to the non-$L^2$ slowly decaying tail \fref{gammapositif}. Let a radial cut off function $\chi_A(r)=\chi\left(\frac{r}{A}\right)$ with $\chi(r)=1$ for $0\leq r\leq 1$ and $\chi(r)=0$ for $r\geq 2$ with the choice:
\be
\label{choixdea}
A=A(t)=e^{2a\frac{\theta(2)}{b(t)}} \ \ \mbox{so that} \ \ \gab^{-\frac{a}{2}}\leq A\leq \gab^{-\frac{3a}{2}},
\ee
for some parameter $a>0$ small enough to be chosen later and which depends on $\eta$. Let $$\zh=\chi_A\zeta_b=\zhre+i\zhim.$$ $\zh$ still satisfies size estimates of Lemma \ref{outgoingradiation} and is moreover in $H^1$ with estimate:
\be
\label{estnormeldeuxzetat}
|(1+|y|)^{10}(|\zh|+|\nabla \zh|)|_{L^2}^2+\left|(1+|y|)^{10}\left(\left|\frac{\partial \zh}{\partial b}\right|+\left|\nabla \frac{\partial \zh}{\partial b}\right|\right)\right|_{L^2}^2\leq \gab^{1-C\eta}.
\ee
From \fref{zeta}, the equation satisfied by $\zt_b$ is now:$$\Delta\zh-\zh+ibD\zh=\Psi^{(0)}_b+F$$ with
\be
\label{eqfzeta}
F=(\Delta \chi_A)\zeta_b+2\nabla\chi_A\cdot\nabla\zeta_b+iby\cdot\nabla\chi_A\zeta_b.
\ee
We then consider the new profile and dispersion: $$\qbh=Q_b+\zh, \ \ \eh=\e-\zh$$ and claim the following refined local virial estimate for $\eh$:

\begin{Lemma}[Refined virial estimate for $\eh$]
\label{proprefinedviriel}
There holds for some universal constant $c_2>0$:
\be
\label{refinedvirial}
\left\{f_1(s)\right\}_s\geq c_2\left(\int |\nabla \eh|^2+\normeldloch+\gab\right)-\frac{1}{c_2}\left(\sigma_c+\int_A^{2A}|\e|^2\right),
\ee
with
\be 
\label{fun}
f_1(s)=-\frac{1}{2}Im\left(\int y\nabla\qbh\overline{\qbh}\right)-Im\left(\hat{\zeta}_b,\overline{\Lambda \hat{Q}_b}\right)+(\e_2,\Lambda \zhre)-(\e_1,\Lambda\zhim).
\ee
\end{Lemma}

{\bf Proof of Lemma \ref{proprefinedviriel}}\\

The proof is parallel to the one of the local virial estimate \fref{lowerboundbs} up to the estimate of the remaining leading order liner term \fref{remaingicrudebound} for which will use a sharp flux computation based on \fref{gammapositif}.\\
  
  {\bf step 1} Estimate on $\hat{\Psi}_b$\\

Let 
\be
\label{dehaptspib}
\hat{\Psi}_b=-i\sigma_c \mu_b\frac{\partial \hat{Q}_b}{\partial b}-\Delta\qbh+\qbh -ib\Lambda \qbh-\qbh|\qbh|^{p-1}=\hat{\Psi}_b^{(2)}-F
\ee with explicitly: 
$$\hat{\Psi}^{(2)}_b =  -i\sigma_c\mu_b\frac{\partial \zh}{\partial b}+ib\sigma_c\zh+\Psi_b^{(1)}-  \left[(\qb+\zh)|\qb+\zh|^{p-1}-Q_b|Q_b|^{p-1}\right],
$$where $\Psi_b^{(1)}$ is given by \fref{defpsib}. Observe from \fref{etglobalepsitwo}, \fref{etglobalepsitwobis}, \fref{estnormeldeuxzetat} and the degeneracy on compact sets \fref{estgammacompact} that:
\be
\label{estrestepsibtwo}
|(1+|y|)^{10}(|\hat{\Psi}_b^{(2)}|+|\nabla \hat{\Psi}_b^{(2)}|)|_{L^2}^2\leq \gab^{1+c}+C\sigma_c^2
\ee
for some universal constants $c,C>0$.\\

{\bf step 2} Rerunning the local virial estimate bound.\\

We now turn to the proof of \fref{refinedvirial} and propose a small short cut with respect to the analysis in \cite{MR4}. Let us indeed rerun exactly step 1 of the proof of Proposition \ref{propviriel} with the new profile $\qbh$ and variable $\eh$. Recall indeed that the whole algebra is completely intrinsic to the equation of the profile and relies only on the definition of respectively  \fref{defpsib}, \fref{dehaptspib}. Moreover, the algebraic formula \fref{firstexpressionbsbisbis} does not rely on the choice of orthogonality conditions for $\e$ --which no longer hold for $\eh$--, and hence:
\bea
\label{cnococ}
\nonumber   & & b_s\left(\left(\dtdbh,\Lambda\Sh\right)-\left(\dsdbh,\Lambda\Th\right)-\left(\eh_2,\Lambda\dsdbh\right)+\left(\e_1,\Lambda \dtdbh\right)\right)+\left\{(\eh_2,\Lambda\Sh)-(\eh_1,\Lambda\Th)\right\}_s=\\
\nonumber  & & \delta(p)\sigma_c+\sigma_c(2E(\qbh)+|\qbh|_{L^2}^2+\tgamma_s|\qbh|_{L^2}^2)-2\lambda^{2(1-\sigma_c)}E_0\\
& - & (\eh_1,Re(\Lambda\hat{\Psi}_b))-(\eh_2,Im(\Lambda\hat{\Psi}_b))+ H_p(\eh,\eh)-\tgamma_s\left\{(\eh_1,\Lambda\Sh)+(\eh_2,\Lambda\Th)\right\}\\
 \nonumber & - & \xsl\cdot\left\{(\eh_2,\nabla\Lambda\Sh)-(\eh_1,\nabla\Lambda \Th)\right\}+ E(\eh,\eh)+(G_1(\eh),\Lambda \Sh)+(G_2(\eh),\Lambda\Th)\\
\nonumber  & - & \left(\lsl+b\right)\left\{(\eh_2,\Lambda^2\Sh+2\sigma_c\Lambda\Sh)-(\eh_1,\Lambda^2\Th+2\sigma_c\Lambda\Th)\right\}-\frac{2}{p+1}\int F(\eh).
\eea
We first observe after an integration by parts that:
$$\left(\dtdbh,\Lambda\Sh\right)-\left(\dsdbh,\Lambda\Th\right)=-\frac{1}{2}\frac{d}{db}Im\left(\int y\cdot\nabla \qbh\overline{\qbh}\right)-\sigma_c \left[\left(\dtdbh,\Sh\right)-\left(\dsdbh,\Th\right)\right].$$ Using the orthogonality relation \fref{orthe4b}, we compute:
\bee
& & -\frac{1}{2}Im\left(\int y\cdot\nabla \qbh\overline{\qbh}\right)+(\eh_2,\Lambda\Sh)-(\eh_1,\Lambda\Th) \\
& = & -\frac{1}{2}Im\left(\int y\cdot\nabla \qbh\overline{\qbh}\right)-2Im\left(\hat{\zeta}_b,\overline{\Lambda \hat{Q}_b}\right)+(\e_2,\Lambda \zhre)-(\e_1,\Lambda\zhim)\\
& = & f_1.
\eee
Next, all the terms in \fref{cnococ} are treated like for the proof of \fref{lowerboundbs} except the linear term involving $\hat{\Psi}_b$. Arguing as in \cite{MR4} and using in particular the degeneracy of $\zeta_b$ on compact sets \fref{estgammacompact} to treat the scalar products terms in $\eh$, we arrive at the following preliminary estimate:
\bea
\label{ceiohcioehf}
\nonumber \{f_1\}_s & \geq & c_0\left(\int |\nabla \eh|^2+\normeldloch\right)-C\lambda^{2(1-\sigma_c)}E_0-C\sigma_c-\gab^{1+z_0}\\
& - & (\eh_1,Re(\Lambda\hat{\Psi}_b))-(\eh_2,Im(\Lambda\hat{\Psi}_b))
\eea
for some universal constants $C,z_0>0$ and with $f_1$ given by \fref{fun}. It remains to estimate the leading order linear term. We first estimate from \fref{dehaptspib}, \fref{estrestepsibtwo}:
\bea
 \label{cndkoneonveo}
 -   (\eh_1,Re(\Lambda\hat{\Psi}_b))-(\eh_2,Im(\Lambda\hat{\Psi}_b)) & \geq &  (\eh_1,Re(DF)+(\eh_2,Im(DF))\\
 \nonumber & - & \frac{c_0}{10}\left(\int |\nabla \eh|^2+\normeldloch\right)-\gab^{1+c_0}-C\sigma_c
\eea
To estimate the remaining linear term, we proceed as in \cite{MR4} and split $\eh=\e-\zh$:
$$ (\eh_1,Re(DF))+(\eh_2,Im(DF))= (\e_1,Re(DF))+(\e_2,Im(DF))- (\zhre,Re(DF))-(\zhim,Im(DF)).$$ The last term is the flux term for which the following lower bound can be derived from \fref{gammapositif}: 
\be
\label{computflux}
-(\zhre,Re(DF))-(\zhim,Im(DF))\geq c_0\gab.
\ee 
The other term is estimated from Cauchy Schwarz and a sharp estimate on $F$ from \fref{eqfzeta} and \fref{gammapositif}, \fref{controlegrad}: 
\bea
\label{vcmlmvdep}
\nonumber |(\e_1,Re(DF))+(\e_2,Im(DF))| & \lesssim & \left(\int_{A}^{2A}|\e|^2\right)^{\frac{1}{2}} \left(\int_{A}^{2A}|F|^2\right)^{\frac{1}{2}}\\
& \lesssim & \frac{c_0}{10}\gab +\frac{10}{c_0}\int_{A}^{2A}|\e|^2.
\eea
We refer to \cite{MR4}, step 4 of the proof of Lemma 6, for a detailed proof of \fref{computflux}, \fref{vcmlmvdep}. Injecting \fref{computflux}, \fref{vcmlmvdep} into \fref{cndkoneonveo} and \fref{ceiohcioehf}  now yields \fref{refinedvirial}.\\
This concludes the proof of Lemma \ref{proprefinedviriel}.


\subsection{Computation of the $L^2$ flux}


We now turn to the computation of $L^2$ fluxes which are the key to get upper bounds on the far away localized $L^2$ term which appears in the RHS of \fref{refinedvirial}. The obtained identity displays new features with respect to the analysis in \cite{MR4} which reflect the $L^2$ super critical nature of the problem. We introduce a radial non negative cut off function $\phi(r)$ such that $\phi(r)=0$ for $r\leq \frac{1}{2}$, $\phi(r)=1$ for $r\geq 3$, $\frac{1}{4}\leq \phi'(r)\leq \frac{1}{2}$ pour $1\leq r\leq 2$, $\phi'(r)\geq 0$. We then let $$\phi_A(s,r)= \phi\left(\frac{r}{A(s)}\right),$$ with $A(s)$ given by (\ref{choixdea}).

\begin{Lemma}[$L^2$ fluxes]
\label{lemmalinearcontrol}
There holds for some universal constant $c_3,z_0>0$ and $s\geq 0$:
\be
\label{estlineaireinfinity}
\frac{1}{\lambda^{2\sigma_c}}\left\{\lambda^{2\sigma_c}\int \phi_A |\e|^2\right\}_s\geq c_3b\int_A^{2A}|\e|^2-\gab^{1+z_0}-\gab^{\frac{a}{2}}\int |\nabla\e|^2.
\ee
\end{Lemma}

{\bf Proof of Lemma \ref{lemmalinearcontrol}}\\

Take a smooth cut off function $\chi(t,x)$. We integrate by parts on \fref{nls} to compute the flux of $L^2$ norm: $$\frac{1}{2}\left\{\int \chi(t,x)|u(t,x)|^2dx\right\}_t=\frac{1}{2}\int \partial_t\chi(t,x)|u(t,x)|^2dx+Im\left(\nabla\chi\cdot \nabla u\overline{u}\right).$$ We apply this with $\chi(t,x)=\phi_A(\frac{x-x(t)}{\lambda(t)})$ and inject the decomposition \fref{decompou}. Recall by construction that $Q_b(y)=\sigma_cT_b(y)$ for $|y|\geq \frac{2}{b}$ which is uniformly exponentially decreasing and hence using \fref{comparaisonbalphaboot}, its contribution near $A$ generates terms which are negligible with respect to the leading order $\gab$. We get after a bit of algebra using also the bootstrap estimates:
\bea
\label{displinearalgebra}
\nonumber & & \frac{1}{2\lambda^{2\sigma_c}}\left\{\lambda^{2\sigma_c}\int \phi_A |\e|^2\right\}_s\geq \frac{b}{2}\int y\cdot\nabla \phi_A |\e|^2+\frac{1}{2}\int \frac{\partial \phi_A}{\partial s} |\e|^2 +  Im\left(\int \nabla\phi_A\cdot\nabla\e\overline{\e}\right)\\
 & - & \frac{1}{2}\left(\lsl+b\right)\int y\cdot\nabla \phi_A |\e|^2-\frac{1}{2}\xsl\cdot\int\nabla \phi_A |\e|^2-\gab^{1+z_0}
\eea
for some universal constant $z_0>0$. From the choice of $\phi$:
\be
\label{contrfluxdeuxun}
10\int \phi'\left(\frac{y}{A}\right) |\e|^2\geq \frac{1}{A}\int y\cdot\nabla \phi\left(\frac{y}{A}\right)|\e|^2\geq \frac{1}{10}\int \phi'\left(\frac{y}{A}\right) |\e|^2\geq \frac{1}{40}\int_{A}^{2A}|\e|^2,
\ee
and also from the choice of $A$ and the control of the geometrical parameters:
\be
\label{vjpovjeopjep}
\left|\left(\lsl+b\right)\int y\cdot\nabla \phi_A |\e|^2\right|+\left|\xsl\cdot\int\nabla \phi_A |\e|^2\right|+\left|\int \frac{\partial \phi_A}{\partial s} |\e|^2\right|\lesssim \frac{b}{1000}\int \phi'\left(\frac{y}{A}\right)|\e|^2.
\ee
Moreover:
\bea
\label{controlfluxun}
\nonumber & & \left| Im\left(\int \nabla\phi_A\cdot\nabla \e\overline{\e}\right)\right|  = \left| Im\left(\int \frac{1}{A}\nabla \phi\left(\frac{y}{A}\right)\cdot\nabla \e\overline{\e}\right)\right|\\
\nonumber & \leq & \frac{1}{A}\left(\int|\nabla\e|^2\right)^{\frac{1}{2}}\left(\int \phi'\left(\frac{y}{A}\right)|\e|^2\right)^{\frac{1}{2}}\\
\nonumber & \leq & \frac{40}{bA^2}\int |\nabla \e|^2+\frac{b}{40}\int \phi'\left(\frac{y}{A}\right)|\e|^2\\
 & \leq & \frac{b}{40}\int \phi'\left(\frac{y}{A}\right)|\e|^2+\gab^{\frac{a}{2}}\int|\nabla \e|^2.
\eea
Injecting \fref{contrfluxdeuxun}, \fref{vjpovjeopjep} and \fref{controlfluxun} into \fref{displinearalgebra} yields \fref{estlineaireinfinity} and concludes the proof of Lemma \ref{lemmalinearcontrol}.


\subsection{$L^2$ conservation law and second monotonicity formula}


We now couple the estimates \fref{refinedvirial} and \fref{estlineaireinfinity} together with the $L^2$ conservation law to derive a new monotonicity formula which completes the dynamical information given by \fref{lowerboundbs}.

\begin{Prop}[Second monotonicity formula]
\label{lemmadiffineqb}
There holds for some universal constant $c_4>0$:
\be
\label{diffineqb}
 -\left\{{\mathcal{J}}\right\}_s\geq c_4b\left(\gab+\int|\nabla \eh|^2+\int|\eh|^2e^{- |y|}\right)-\frac{b}{c_4}\sigma_c,
\ee
with
\bea
\label{calftrois}
{\mathcal{J}}(s) & = & \left(\int |\qb|^2-\int Q_p^2\right)+2(\e_1,\S)+2(\e_2,\T)+\int (1-\phi_A) |\e|^2\\
 \nonumber& - & c_3c_2\left( b\tilde{f}_1(b)-\int_0^b\tilde{f}_1(v)dv+b\{(\e_2,\Lambda \zhre)-(\e_1,\Lambda\zhim)\}\right),
\eea
where $c_3,c_2$ are the universal small constants involved in \fref{refinedvirial}, \fref{estlineaireinfinity}, and:
\be
\label{eqfuntilde}
\tilde{f}_1(b)=\frac{1}{2}Im\left(\int y\nabla\hat{Q_b}\overline{\hat{Q_b}}\right)+(\T,\Lambda\Sh)-(\S,\Lambda \Th).
\ee
\end{Prop}

{\bf Proof of Proposition \ref{lemmadiffineqb}}\\

{\bf step 1} Coupling \fref{refinedvirial} and \fref{estlineaireinfinity}.\\

Let us multiply \fref{refinedvirial} by $bc_2$:
\be
\label{cvnevnevnepp}
bc_2^2\left(\int |\nabla \eh|^2+\int |\eh|^2e^{-|y|}+\gab\right)\leq b\{c_2f_1\}_s+b\int_A^{2A}|\e|^2+b\sigma_c.
\ee
We then integrate by parts in time using \fref{fun}:
\bee
 b\{c_2f_1\}_s & = & c_2\left\{ b\tilde{f}_1(b)-\int_0^b\tilde{f}_1(v)dv+b(\e_2,\Lambda \zhre)-b(\e_1,\Lambda\zhim)\right\}_s\\
 & - & c_2b_s\left\{(\e_2,\Lambda \zhre)-(\e_1,\Lambda\zhim)\right\}
 \eee
and estimate from \fref{smalnesscoeffpairb} and \eqref{estnormeldeuxzetat}:
\bea
&&\nonumber\left|c_2b_s\left\{(\e_2,\Lambda \zhre)-(\e_1,\Lambda\zhim)\right\}\right|\\
\nonumber &\leq &\left|c_2b_s\left\{(\eh_2,\Lambda \zhre)-(\eh_1,\Lambda\zhim)\right\}\right|+\left|c_2b_s\left\{(\zhre,\Lambda \zhim)-(\zhim,\Lambda\zhre)\right\}\right|\\
& \leq &\gab^{1+z_0}+\frac{bc_2^2}{10}\left(\int |\nabla \eh|^2+\int |\eh|^2e^{-|y|}+\gab\right).
\eea 
Injecting this into \fref{cvnevnevnepp} yields:
\bee
\frac{bc_2^2}{2}\left(\int |\nabla \eh|^2+\int |\eh|^2e^{-|y|}+\gab\right) & \leq &  c_2\left\{ b\tilde{f}_1(b)-\int_0^b\tilde{f}_1(v)dv+b(\e_2,\Lambda \zhre)-b(\e_1,\Lambda\zhim)\right\}_s\\
& + & b\int_A^{2A}|\e|^2+b\sigma_c.
\eee
We now inject the control of $L^2$ fluxes \fref{estlineaireinfinity} and obtain for $a>C\eta$:
\bea
\label{cnovneoeuyo}
\nonumber \frac{bc_3c_2^2}{4}\left(\int |\nabla \eh|^2+\int |\eh|^2e^{-|y|}+\gab\right) & \leq &  c_3c_2\left\{ b\tilde{f}_1(b)-\int_0^b\tilde{f}_1(v)dv+b(\e_2,\Lambda \zhim)-b(\e_1,\Lambda\zhre)\right\}_s\\
& + & \frac{1}{\lambda^{2\sigma_c}}\left\{\lambda^{2\sigma_c}\int \phi_A |\e|^2\right\}_s+bc_3\sigma_c.
\eea

{\bf step 2} Injection of the $L^2$ conservation law.\\

We now rewrite the $L^2$ conservation law as follows:
$$\int|u_0|^2=\lambda^{2\sigma_c}\left(\int|Q_b|^2+2Re(\e,\overline{Q_b})+\int|\e|^2\right),$$
which yields:
\bee
& & \frac{1}{\lambda^{2\sigma_c}}\left\{\lambda^{2\sigma_c}\int \phi_A |\e|^2\right\}_s\\ 
& = & -\frac{1}{\lambda^{2\sigma_c}}\left\{\lambda^{2\sigma_c}\left[\int (1-\phi_A) |\e|^2+2Re(\e,\overline{Q_b})+\int|Q_b|^2\right]\right\}_s\\
& = & -\left\{\int (1-\phi_A) |\e|^2+2Re(\e,\overline{Q_b})+\int|Q_b|^2\right\}_s\\
& & -2\sigma_c\lsl\left[\int (1-\phi_A) |\e|^2+2Re(\e,\overline{Q_b})+\int|Q_b|^2\right].
\eee
Now the Hardy type bound
\be
\label{estintermjfjls}
\int (1-\phi_A) |\e|^2\leq CA^3\left(\int |\nabla \e|^2+\normeldloc\right)
\ee
together with the choice of $A$ (\ref{choixdea}), the bootstrap bound \fref{estenboot} and the control \fref{smalnesscoeffpairb} yield the rough bound: $$\left|\sigma_c\lsl\left[\int (1-\phi_A) |\e|^2+2Re(\e,\overline{Q_b})+\int|Q_b|^2\right]\right|\lesssim Cb\sigma_c.$$ 
We may thus rewrite \fref{cnovneoeuyo} as:
\bee
& & \frac{bc_3c_2^2}{4}\left(\int |\nabla \eh|^2+\int |\eh|^2e^{-|y|}+\gab\right)\leq Cb\sigma_c\\
& +&   \left\{ -\left(\int|Q_b|^2-\int Q_p^2\right)-2Re(\e,\overline{Q_b})-\int (1-\phi_A) |\e|^2\right .\\
& + & \left . c_3c_2\left[b\tilde{f}_1(b)-\int_0^b\tilde{f}_1(v)dv+b(\e_2,\Lambda \zhim)-b(\e_1,\Lambda\zhre)\right]    \right\}_s,
\eee
which is \fref{diffineqb}. This concludes the proof of Proposition \ref{lemmadiffineqb}.


\section{Existence and stability of the self similar regime}
\label{sectionproofthm}


This section is devoted to the proof of the main Theorem \ref{thmmain}. We first show how the {\it coupling} of the monotonicity formulae \fref{lowerboundbs}, \fref{diffineqb} implies a dynamical trapping of $b$ and a uniform bound on $\e$ which allows us to close the bootstrap Proposition \ref{propboot}. We then conclude the proof of Theorem \ref{thmmain} as a simple consequence of these uniform bounds.


\subsection{Closing the bootstrap}
\label{sectionclosingboot}


We are now in position to close the bootstrap and conclude the proof of Proposition \ref{propboot}.\\

{\bf Proof of Proposition \ref{propboot}}\\

{\bf step 1} Pointwise bound on $\e$.\\

Let us start with the proof of the pointwise bound on $\e$ \fref{estebootfinal}. We argue by contradiction and assume that there exists $s_2\in [s_0,s^*]$ such that: $$ \int|\nabla \e(s_2)|^2+\int|\e(s_2)|^2e^{-|y|}> \Gamma_{b(s_2)}^{1-9\nu_0}.$$ A simple continuity argument based on the initialization of the bootstrap estimate \fref{esteinit} implies that there exists $[s_3,s_4]\subset[s_0,s^*]$ such that: 
\be
\label{eststrhreesfour}
 \int|\nabla \e(s_3)|^2+\int|\e(s_3)|^2e^{-|y|}= \Gamma_{b(s_3)}^{1-7\nu_0}, \ \  \int|\nabla \e(s_4)|^2+\int|\e(s_4)|^2e^{-|y|}= \Gamma_{b(s_4)}^{1-9\nu_0},
  \ee and
 \be
 \label{vodvodovintervalsthreessouf}
 \forall s\in [s_3,s_4], \ \  \int|\nabla \e(s)|^2+\int|\e(s)|^2e^{-|y|}\geq \Gamma_{b(s)}^{1-7\nu_0}.
 \ee
 From \fref{vodvodovintervalsthreessouf} and the first virial monotonicity \fref{lowerboundbs}, we have: $\forall s\in [s_3,s_4]$, $$b_s\geq c_1( \Gamma_{b}^{1-7\nu_0}-\gab^{1-\nu^2_0})>0$$ and hence 
 \be
 \label{monotonicityb}
 b(s_4)\geq b(s_3).
 \ee 
 On the other hand, using the lower bound $$\int|\nabla \eh|^2+\int|\eh|^2e^{-|y|}\geq \frac{1}{2}\left(\int|\nabla \e|^2+\int|\e|^2e^{-|y|}\right)-\gab^{1-\nu_0}$$ together with \fref{vodvodovintervalsthreessouf}, \fref{comparaisonbalphaboot} and the second monotonicity formula \fref{diffineqb}, there holds: $\forall s\in [s_3,s_4]$, 
 $$-\mathcal J_s\geq \frac{bc_4}{2}(\Gamma_{b}^{1-7\nu_0}-\gab^{1-6\nu_0})\geq 0$$ and hence
 \be
 \label{monotninicityj}
 \mathcal J(s_4)\leq \mathcal J(s_3).
 \ee
We now claim the following upper and lower control of $\mathcal J$: 
\be
\label{contrfdeux}
{\mathcal{J}}(s)-f_2(b(s)) \ \  \left  \{ \begin{array}{ll}
        \geq -\gab^{1-Ca}+\frac{1}{C}\left(\int |\nabla\e|^2+\normeldloc\right)-C\sigma_c,\\
         \leq  CA^3\left(\int |\nabla\e|^2+\normeldloc\right)+\gab^{1-Ca}+C\sigma_c,
         \end{array}
\right .
\ee
where $f_2$ given by 
\be
\label{deffrdexu}
f_2(b,\sigma_c)=\left(\int |\qb|^2-\int Q_p^2\right)-c_3c_2\left( b\tilde{f}_1(b)-\int_0^b\tilde{f}_1(v)dv\right) 
\ee
satisfies
\be
\label{lowergfdeux}
\forall b^*>b_2>b_1, \ \ f_2(b_2)\geq f_2(b_1)-C\sigma_c, \ \ \frac{1}{C}b_1^2-C\sigma_c\leq f_2(b_1)\leq C(b_1^2+\sigma_c).
\ee
Let us assume \fref{contrfdeux}, \fref{lowergfdeux}. Then  \fref{monotninicityj}, \fref{contrfdeux} imply:
\bee
& & f_2(b(s_4))-\Gamma_{b(s_4)}^{1-Ca}+\frac{1}{C}\left(\int |\nabla\e(s_4)|^2+\int|\e(s_4)|^2e^{-|y|}\right)-C\sigma_c\leq  \mathcal J(s_4)\leq \mathcal J(s_3)\\
& \leq & f_2(b(s_3))+CA^3(s_3)\left(\int |\nabla\e(s_3)|^2+\int|\e(s_3)|^2e^{-|y|}\right)+\Gamma_{b(s_3)}^{1-Ca}+C\sigma_c
\eee
and hence from the monotonicity \fref{monotonicityb}, \fref{lowergfdeux} and the controls \fref{eststrhreesfour}, \fref{comparaisonbalphaboot}:
$$
\Gamma_{b(s_4)}^{1-9\nu_0} =\int|\nabla \e(s_4)|^2+\int|\e(s_4)|^2e^{-|y|}
 \leq   C\Gamma_{b(s_4)}^{1-Ca}+C\Gamma_{b(s_3)}^{1-7\nu_0-Ca}+C\sigma_c\leq  \Gamma_{b(s_4)}^{1-8\nu_0}
$$
for $a=C\eta>0$ and $\eta$ chosen small enough, a contradiction which concludes the proof of \fref{estebootfinal}.\\
{\it Proof of \fref{contrfdeux}, \fref{lowergfdeux}}: It is a standard consequence of the coercivity of the linearized energy with our choice of orthogonality conditions, \cite{MR4}. Indeed, we rewrite $\mathcal J$ given by \fref{calftrois} using the conservation of energy \fref{energyconservation} and the orthogonality condition \eqref{orthe4b}:
\bee
\mathcal J & = &f_2(b,\sigma_c)+2\left(\e_1,Re(\Psi_b)-\sigma_c\mu_b\dtdb\right)+2\left(\e_2,Im (\Psi_b)+\sigma_c\mu_b\dsdb\right) \\
& - & c_3c_2b\{(\e_2,\Lambda \zhre)-(\e_1,\Lambda\zhim)\}+ 2E(Q_b)-2\lambda^{2(1-\sigma_c)}E_0\\
& + & (M_+(\e),\e_1)+(M_-(\e),\e_2)-\int\phi_A|\e|^2-\frac{2}{p+1}\int F(\e).
\eee
The upper bound in \fref{contrfdeux} now follows from the Hardy bound \fref{estintermjfjls} and the degeneracy \fref{calculhamiltonian}. For the lower bound, we recall the following coercivity of the linearized energy which holds true for $A$ large enough i.e. $|b|\leq b^*$ small enough: 
\bee
& & (M_+( \e),\e_1)+(M_-(\e),\e_2)-\int\phi_A|\e|^2  \geq c_3\left(\int |\nabla\e|^2+\normeldloc\right)\\
\nonumber & - & \frac{1}{c_3}\left\{(\e_1,Q_{p_c})^2+(\e_1,|y|^2Q_{p_c})^2+(\e_1,yQ_{p_c})^2+(\e_2,D^2Q_{p_c})^2\right\},
\eee
for some universal constant $c_3>0$, see  Appendix D in \cite{MR4}. The choice of orthogonality conditions together with the degeneracy \eqref{controlelenergyb} now yield \fref{contrfdeux}. \fref{lowergfdeux} is now a direct consequence of \fref{qbtsupcritbis}, \fref{propmb}. Indeed, first observe from \fref{eqfuntilde}, \fref{degenviriel} and the smallness of the radiation given by Lemma \ref{outgoingradiation} that:
$$\tilde{f_1}(b)=M_1(b)+O(\sigma_c) \ \ \mbox{with} \ \ M_1(0)=0 \ \ \mbox{and} \ \ \left|\frac{d M_1(b)}{db}\right|\lesssim C(p) $$ for some universal constant $C(p)>0$. Hence \fref{qbtsupcritbis}, \fref{propmb} imply:
\bee
f_2(b,\sigma_c)& = & \left(\int |\qb|^2-\int Q_p^2\right)-c_3c_2\left( b\tilde{f}_1(b)-\int_0^b\tilde{f}_1(v)dv\right)\\
& = & \tilde{M}(b)+O(\sigma_c)
\eee with $$\frac{d\tilde{M}(b)}{db}\geq b\left(c_0(p)-Cc_3c_2\right)\geq \frac{b}{2}c_0(p_c)$$ provided the constants $c_2$, $c_3$ in \fref{refinedvirial}, \fref{estlineaireinfinity} have been chosen small enough, and the monotonicity \fref{lowergfdeux}  follows.\\

{\bf step 2} Dynamical trapping of $b$.\\

We now turn to the core of the argument which is the dynamical trapping of $b$ \fref{comparaisonbalphabootalpha}. We recall from \fref{gammapositif} that $$e^{-\frac{2\theta(2)}{b}(1+C\eta)}\leq \Gamma_b\leq e^{-\frac{2\theta(2)}{b}(1-C\eta)}.$$
Assume that there exists $s_5\in [s_0,s^*]$ such that $\sigma_c\geq \Gamma_{b(s_5)}^{1-\nu^4_0}$, then from \fref{comparaisonbalphainit} and a simple continuity argument, consider $s_6\in[s_0,s_5)$ the largest  time such that $\sigma_c=\left(e^{-2\frac{\theta(2)}{b(s_6)}}\right)^{1-\nu^5_0}$, then $b_s(s_6)\leq 0$ by construction while from \fref{lowerboundbs}:
$$b_s(s_6)\geq c_1\left(\sigma_c-\Gamma_{b(s_6)}^{1-\nu^6_0}\right)\geq c_1\left(\left(e^{-2\frac{\theta(2)}{b(s_6)}}\right)^{1-\nu^5_0}-\Gamma_{b(s_6)}^{1-\nu^6_0}\right)>0$$ and a contradiction follows.\\
Similarly, assume that there exists $s_7\in [s_0,s^*]$ such that 
\be
\label{neoeoheo}
\sigma_c\leq \Gamma_{b(s_7)}^{1+\nu^4_0}.
\ee From \fref{calftrois}, the pointwise bounds \fref{comparaisonbalphaboot}, \fref{estenboot}, and the value of the $L^2$ norm of $Q_b$ given by \fref{qbtsupcritbis}, \fref{propmb}, there holds: 
\be
\label{estleucp}
|\mathcal J-d_0b^2|\lesssim \nu_0^{100}b^2
\ee for some universal constant $d_0>0$. Hence \fref{comparaisonbalphainit}, \fref{neoeoheo} imply: $$\sigma_c\geq \left(\Gamma_{\sqrt{\frac{\mathcal J(s_0)}{d_0}}}\right)^{1+\nu^9_0}, \ \  \sigma_c\leq \left(\Gamma_{\sqrt{\frac{\mathcal J(s_7)}{d_0}}}\right)^{1+\nu^5_0}.$$ We then consider -using again \fref{gammapositif}- the largest time $s_8\in[s_0,s_7]$ such that $\sigma_c=\left(e^{-2\frac{\theta(2)}{\sqrt{\frac{\mathcal J(s_8)}{d_0}}}}\right)^{1+\nu^6_0}$, then $(\mathcal J)_s(s_8)\geq 0$ by definition while from \fref{diffineqb}, \fref{estleucp}, \fref{gammapositif}:
$$ -\left\{{\mathcal{J}}\right\}_s(s_8)\geq b(s_8)\left[c_4\left(\Gamma_{\sqrt{\frac{\mathcal J(s_8)}{d_0}}}\right)^{1+\nu^8_0}-\frac{1}{c_4}\left(e^{-2\frac{\theta(2)}{\sqrt{\frac{\mathcal J(s_8)}{d_0}}}}\right)^{1+\nu^6_0}\right]>0,$$ and a contradiction follows. This concludes the proof of \fref{comparaisonbalphabootalpha}.\\

{\bf step 3} Control of the scaling parameter.\\

We now turn to the control of the scaling parameter $\lambda(t)$. From \fref{comparaisonbalphainit}, the dynamical trapping of $b$ \fref{comparaisonbalphabootalpha} implies: 
\be
\label{uniformcontorl}
\forall t\in [0,T^*), \ \ (1-\nu^3_0)b_0\leq b(t)\leq (1+\nu^3_0)b_0.
\ee
Hence the upper bound \fref{smalnesscoeffpairb}, the control \fref{estebootfinal} and \fref{uniformcontorl} ensure: 
\be
\label{pointwisevounfl}
\forall t\in[0,T^*), \ \ 0<(1-2\nu^3_0)b_0\leq-\lsl=-\lambda_t\lambda\leq (1+2\nu^3_0)b_0.
\ee In particular, $\lambda$ is nonincreasing while $b(t)$ is trapped from \fref{uniformcontorl} and thus \fref{lambdasmallimit}, \fref{initienergyomentum} imply \fref{estlambdabootfinal}, \fref{energyomentumfinal}.\\

{\bf step 4} Control of the solution in $\dot{H}^{\sigma}$.\\

It remains to close the $\dot{H}^{\sigma}$ estimate \fref{esteisigmafinal} which is the key to the control of the nonlinear term \fref{estnonlinearterm}. We use here the fact that the blow up is self similar and strictly $H^1$ subcritical so that $$\sigma_p=\frac{N}{2}-\frac{N}{p+1}>\sigma_c=\frac{N}{2}-\frac{2}{p-1}.$$ In other words, norms above scaling can be controlled dynamically in the bootstrap as was for example observed in \cite{RS}, \cite{RodSter}, \cite{RRS}.\\
It is more convenient here to work in original variables. Consider the decomposition 
$$u(t,x)=Q_{sing}(t,x)+\tilde{u}(t,x)=\frac{1}{\lambda^{\frac{2}{p-1}}(t)}\left(Q_b+\e\right)\left(t,\frac{x-x(t)}{\lambda(t)}\right)e^{i\gamma(t)},
$$
then first observe by rescaling and the trapping of $b$ \fref{uniformcontorl} that \fref{esteisigmafinal} is implied by:
\be
\label{tobeprovedhsigma}
|\tilde{u}|^2_{\dot{H}^{\sigma}}\leq \frac{\Gamma_{b}^{1-45\nu_0}}{\lambda^{2(\sigma-\sigma_c)}}.
\ee
To prove \fref{tobeprovedhsigma}, we write down the equation for $\tilde{u}$ and use standard Strichartz estimates, see \cite{Cazenavebook}, for the linear Schr\"odinger flow. Indeed, the equation for $\tilde{u}$ is: 
$$i\partial_t\tilde{u}+\Delta \tilde{u}=-\mathcal E-f(\tilde{u}) $$ with: 
\bee
\mathcal E & = & i\partial_t Q_{sing}+\Delta Q_{sing}+Q_{sing}|Q_{sing}|^{p-1}\\
& = & \frac{1}{\lambda^{2+\frac{2}{p-1}}}\left[ib_s\frac{\partial Q_b}{\partial b}+\Delta Q_b-Q_b+Q_b|Q_b|^{p-1}-i\lsl\Lambda Q_b-i\xsl\cdot\nabla Q_b-\tgamma_sQ_b\right]\left(t,\frac{x-x(t)}{\lambda(t)}\right)e^{i\gamma(t)}\\
& = & \frac{1}{\lambda^{2+\frac{2}{p-1}}}\left[-\Psi_b+i(b_s-\sigma_c\mu_b)\frac{\partial Q_b}{\partial b}-i\left(\lsl+b\right)\Lambda Q_b-i\xsl\cdot\nabla Q_b-\tgamma_sQ_b\right]\left(t,\frac{x-x(t)}{\lambda(t)}\right)e^{i\gamma(t)},
\eee 
and 
\be
\label{deftuilde}
f(\tilde{u})=(Q_{sing}+\tilde{u})|Q_{sing}+\tilde{u}|^{p-1}-Q_{sing}|Q_{sing}|^{p-1}.
\ee
Let $t\in[0,T^*)$, we write down the Duhamel formula on $[0,t]$. Following \cite{Cazenavebook}, we consider the Strichartz pair:
\be
\label{strichartzparit}
r=\frac{N(p+1)}{N+\sigma(p-1)}, \ \ \gamma=\frac{4(p+1)}{(p-1)(N-2\sigma)}, \ \ \frac{2}{\gamma}=\frac{N}{2}-\frac{N}{r}
\ee
 and estimate from Strichartz estimates:
\be
\label{vnognogoeghoe}
 ||\nabla|^{\sigma}\tilde{u}|_{L^{\infty}_{[0,t]}L^2_x} \lesssim   ||\nabla|^{\sigma}\tilde{u}_0|_{L^2}+||\nabla|^{\sigma}\mathcal E|_{L^1_{[0,t]}L^2_x}+||\nabla|^{\sigma}f(\tilde{u})|_{L^{\gamma'}_{[0,t]}L^{r'}_x}.
\ee
Using the fact that $\l$ is nonincreasing together with \eqref{esteinit}, we obtain:
\be\label{est:errinitial}
||\nabla|^{\sigma}\tilde{u}_0|_{L^2}\lesssim \frac{\Gamma_{b_0}^{\frac{1}{2}(1-\nu_0)}}{[\l(t)]^{(\sigma-\sigma_c)}}.
\ee
We claim:
\be
\label{ccnoehoe}
||\nabla|^{\sigma}\mathcal E|_{L^1_{[0,t]}L^2_x}\leq \frac{\Gamma_{b_0}^{\frac{1}{2}(1-15\nu_0)}}{[\lambda(t)]^{(\sigma-\sigma_c)}},
\ee
\be
\label{cnkoenceoneoc}
||\nabla|^{\sigma}f(\tilde{u}(t))|_{L^2}\leq \frac{\Gamma_{b_0}^{\frac{1}{2}(1-41\nu_0)}}{[\lambda(t)]^{(\sigma-\sigma_c)}}
\ee
which together with \fref{uniformcontorl}, \fref{vnognogoeghoe}, \eqref{est:errinitial} yield \fref{tobeprovedhsigma}.\\
{\it Proof of \fref{ccnoehoe}}: From the estimates on the geometrical parameters \fref{smalnesscoeffpairb}, \fref{smallnesscoeffimpairb}, the degeneracy estimates \fref{etglobalepsitwo}, \fref{etglobalepsitwobis}, the pointwise bounds \fref{estebootfinal} and \fref{comparaisonbalphabootalpha}, and the freezing of $b$ \fref{uniformcontorl}, there holds: $\forall t\in [0,t^*]$, 
\bea
\label{esthsignaone}
\nonumber ||\nabla|^{\sigma}\mathcal E(t)|_{L^2} &\lesssim &  \frac{1}{[\lambda(t)]^{2+(\sigma-\sigma_c)}}\left\|-\Psi_b+i(b_s-\sigma_c\mu_b)\frac{\partial Q_b}{\partial b}-i\left(\lsl+b\right)\Lambda Q_b-i\xsl\cdot\nabla Q_b-\tgamma_sQ_b\right\|_{H^1}\\
\nonumber & \leq&  \frac{1}{[\lambda(t)]^{2+(\sigma-\sigma_c)}}\left(\int |\nabla\e|^2+\normeldloc+\gab^{1-11\nu_0}\right)^{\frac{1}{2}}\\
& \leq & \frac{\Gamma_{b_0}^{\frac{1}{2}(1-12\nu_0)}}{[\lambda(t)]^{2+(\sigma-\sigma_c)}}.
\eea
We now observe from the self similar blow up speed estimate \fref{pointwisevounfl}: $\forall q>2$, 
\be
\label{estcjcvjs}
\int_0^t\frac{d\tau}{[\lambda(\tau)]^q}\leq \frac{C}{b_0}\int_{0}^t-\frac{\lambda_t}{[\lambda(\tau)]^{q-1}}d\tau\leq \frac{C}{(q-2)b_0[\lambda(t)]^{q-2}}.
\ee
Integrating \fref{esthsignaone} in time and using \fref{estcjcvjs} yields \fref{ccnoehoe}.\\
{\it Proof of \fref{cnkoenceoneoc}}: This estimate follows  {\it in the bootstrap} using the fact that the blow up is self similar and that $\tilde{u}$ is small in $\dot{H}^{\sigma}$ for $\sigma>\sigma_c$ after renormalization. Indeed, we first claim from standard nonlinear estimates in Besov spaces:
\be
\label{cnionofiorg}
||\nabla|^{\sigma} f(\ut)|_{L^{r'}} \lesssim \frac{1}{\lambda^{p(\tilde{\sigma}-\sigma_c)}}||\nabla|^{\tilde{\sigma}}\e|_{L^2},
\ee
where $\tilde{\sigma}$ is defined by:
\be
\label{stilde}
\tilde{\sigma}=\sigma+\frac{N}{2}-\frac{N}{r}=\sigma+\frac{2}{\gamma}.
\ee
The proof of the estimate \eqref{cnionofiorg} is postponed to the appendix. From direct check, $\sigma<\tilde{\sigma}<1$ providing $\sigma$ has been chosen close enough to the critical scaling exponent $\sigma_c$ himself close enough to 0. We may thus interpolate between $\sigma$ and 1 and use \fref{esteisigmaboot}, \fref{estebootfinal} and \eqref{stilde} to estimate: 
\be
\label{corec}
||\nabla|^{\tilde{\sigma}}\e|_{L^2}\lesssim ||\nabla|^{\sigma}\e|_{L^2}^{\frac{1-\tilde{\sigma}}{1-\sigma}}|\nabla\e|^{\frac{\tilde{\sigma}-\sigma}{1-\sigma}}_{L^2} \lesssim \gab^{\frac{1}{2}\left(1-\left(50-\frac{80}{\gamma(1-\sigma)}\right)\nu_0\right)}.
\ee
Provided $p$ is chosen close enough to $p_c$, and $\sigma$ is chosen close enough to 0, we obtain from and \fref{uniformcontorl}, \eqref{corec}:
\be\label{corec1}
||\nabla|^{\tilde{\sigma}}\e|_{L^2} \lesssim \Gamma_{b_0}^{\frac{1}{2}\left(1-40\nu_0\right)}.
\ee
Injecting this into \fref{cnionofiorg} yields:
\be
\label{cjljcoi}
||\nabla|^{\sigma}f(\ut)|_{L^{\gamma'}_{[0,t]}L^{r'}_x}\lesssim \Gamma_{b_0}^{\frac{1}{2}\left(1-40\nu_0\right)}\left(\int_0^t\frac{d\tau}{[\lambda(\tau)]^{(\tilde{\sigma}-\sigma_c)p\gamma'}}\right)^{\frac{1}{\gamma'}}.
\ee Now from \fref{strichartzparit}, \fref{stilde}, there holds: 
\bea
\label{algebranombre}
\nonumber (\tilde{\sigma}-\sigma_c)p-\frac{2}{\gamma'}& = & p(\sigma-\sigma_c)+\frac{2p}{\gamma}-2\left(1-\frac{1}{\gamma}\right)=p(\sigma-\sigma_c)-2+2\frac{p+1}{\gamma}\\
\nonumber & = & (\sigma-\sigma_c)+(p-1)\left[\sigma-\sigma_c+\frac{N-2\sigma}{2}\right]-2\\
& = & \sigma-\sigma_c.
\eea
In particular, $(\tilde{\sigma}-\sigma_c)p\gamma'=2+\gamma'(\sigma-\sigma_c)>2$, and we may thus inject \fref{estcjcvjs} into \fref{cjljcoi} to conclude:
\be
\label{cnkeobvohoer}
||\nabla|^{\sigma}f(\ut)|_{L^{\gamma'}_{[0,t]}L^{r'}_x}\leq \frac{\Gamma_{b_0}^{\frac{1}{2}(1-41\nu_0)}}{[\lambda(t)]^{\sigma-\sigma_c}}
\ee
for $\nu_0>0$ small enough thanks to $p>1$, this is \fref{cnkoenceoneoc}.\\

This concludes the proof of the bootstrap Proposition \ref{propboot}.


\subsection{Proof of Theorem \ref{thmmain}}


We are now in position to prove Theorem \ref{thmmain}.\\

{\bf Proof of Theorem \ref{thmmain}}\\

Pick $\nu_0>0$, $p\in (p_c,p^*(\nu_0))$ and $u_0\in \mathcal O$ corresponding to $b_0=b^*(p)$ as given by Definition \ref{defdefp}. Note that \fref{loisigma} follows from \fref{comparaisonbalphainit}. Let $u(t)$ be the corresponding solution to \fref{nls} with maximum lifetime interval on the right $[0,T)$, then from Proposition \ref{propboot}, $u(t)$ admits on $[0,T)$ a geometrical decomposition $$u(t,x)=\frac{1}{\lambda^{\frac{2}{p-1}}(t)}(Q_{b(t)}+\e)\left(t,\frac{x-x(t)}{\lambda(t)}\right)e^{i\gamma(t)}$$ which satisfies the estimates of Proposition \ref{propboot}. This implies in particular \fref{uniformbounde}.\\

{\bf step 1} Finite time blow and self similar blow up speed.\\

Recall \fref{pointwisevounfl}: $$\forall t\in [0,T), \ \ (1-\nu^2_0)b_0\leq-\lambda_t\lambda\leq (1+\nu^2_0)b_0.$$ Integrating this in time first from $0$ to $t$ yields: $$\forall t\in [0,T), \ \ (1-\nu^2_0)b_0 t\leq \frac{1}{2}\lambda^2_0 \ \ \mbox{and hence} \ \ T\leq \frac{\lambda^2_0}{2b_0(1-\nu^2_0)}$$ so that the solution blows up in finite time. From the $H^1$ Cauchy theory, $|\nabla u(t)|_{L^2}\to +\infty$ as $t\to T$ and hence from \fref{estebootfinal}, $\lambda(t)\to 0$ as $t\to T$. We thus integrate \fref{pointwisevounfl} from $t$ to $T$ to get: 
$$\forall t\in [0,T], \ \ (1-\nu^2_0)b_0(T-t)\leq \frac{\lambda^2(t)}{2}\leq (1+\nu^2_0)b_0(T-t)$$ which implies \fref{selfsmilar}.\\

{\bf step 2} Convergence of the blow up point.\\

From \fref{smallnesscoeffimpairb} and Proposition \ref{propboot}, we have the rough bound: $$|x_t|=\frac{1}{\lambda}\left|\xsl\right|\lesssim \frac{\Gamma_{b_0}^{\frac{1}{4}}}{\lambda}$$ and thus from \fref{selfsmilar}: $$\int_0^T|x_t|dt\lesssim \int_0^T\frac{\Gamma_{b_0}^{\frac{1}{4}}}{\sqrt{b_0(T-t)}}<+\infty,$$ and \fref{convblowuppiint} follows. We moreover get the convergence rate:
\be
\label{covergencerate}
\left|\frac{x(t)-x(T)}{\lambda(t)}\right|\lesssim \frac{\Gamma_{b_0}^{\frac{1}{4}}}{\sqrt{2b_0(T-t)}}\int_t^T\frac{d\tau}{\sqrt{2b_0(T-\tau)}}\lesssim \Gamma_{b_0}^{\frac{1}{8}}.
\ee

\bs
\ni
{\bf step 3} Strong convergence in $H^{s}$ for $0\leq s<\sigma_c$.\\

We now turn to the proof of \fref{convustarhsigma}. Pick $0\leq s<\sigma_c$. Let $0<\tau\ll 1$ and $0<t<T-\tau$, let $u_{\tau}(t)=u(t+\tau)$ and $v(t)=u_{\tau}(t)-u(t)$, then $v$ satisfies: 
\be
\label{vnonveonoe}
iv_t+\Delta v=u|u|^{p-1}-u_{\tau}|u_{\tau}|^{p-1}.
\ee Consider the Strichartz pair 
$$r=\frac{N(p+1)}{N+s(p-1)}, \ \ \gamma=\frac{4(p+1)}{(p-1)(N-2s)}, \ \ \frac{2}{\gamma}=\frac{N}{2}-\frac{N}{r},$$ then from standard nonlinear estimates in Sobolev spaces --\cite{Cazenavebook}--, we have: 
$$\left||\nabla|^{s}\left[u|u|^{p-1}-u_{\tau}|u_{\tau}|^{p-1}\right]\right|_{L^{r'}}\lesssim ||\nabla|^su|_{L^r}^p+||\nabla|^su_{\tau}|_{L^r}^p\lesssim ||\nabla|^{\tilde{\sigma}}u|_{L^2}^p+||\nabla|^{\tilde{\sigma}}u_{\tau}|_{L^2}^p$$ with $$\tilde{\sigma}=s+\frac{N}{2}-\frac{N}{r}=s+\frac{2}{\gamma}.$$
Now observe that $\tilde{\sigma}\to \frac{N}{N+2}>0$ as $p\to p_c$ and hence $\sigma_c<\sigma<\tilde{\sigma}<1$ from \fref{cheioheoh} for $p$ close enough to $p_c$. We thus estimate from the geometrical decomposition \fref{decompou} and the bounds \fref{estebootfinal}, \fref{esteisigmafinal}:
$$||\nabla|^{\tilde{\sigma}}u|_{L^2}\lesssim \frac{1}{\lambda^{\tilde{\sigma}-\sigma_c}}||\nabla|^{\tilde{\sigma}}(Q_b+\e)|_{L^2}\lesssim \frac{1}{\lambda^{\tilde{\sigma}-\sigma_c}}.$$ We conclude using \fref{algebranombre}:
\bea
\label{nvbnveoeoe}
\nonumber ||\nabla|^{s}\left[u|u|^{p-1}-u_{\tau}|u_{\tau}|^{p-1}\right]|_{L^{\gamma'}_{[t,T-\tau)}L^{r'}} & \lesssim & \left(\int_t^T ||\nabla|^{\tilde{\sigma}}u|_{L^2}^{p\gamma'}\right)^{\frac{1}{\gamma'}}\lesssim \left(\int_t^T \frac{d\tau}{[\lambda(\tau)]^{p\gamma'(\tilde{\sigma}-\sigma_c)}}\right)^{\frac{1}{\gamma'}}\\
& \lesssim & \left(\int_t^T \frac{d\tau}{[\lambda(\tau)]^{2+\gamma'(s-\sigma_c)}}\right)^{\frac{1}{\gamma'}}\to 0 \ \ \mbox{as} \ \ t\to T
\eea
from the scaling law $\lambda(t)\sim \sqrt{2b_0(T-t)}$ and $s-\sigma_c<0.$ By running the standard Strichartz estimates -\cite{Cazenavebook}- on \fref{vnonveonoe}, we conclude that: $$||\nabla|^{s}v|_{L^{\infty}_{[t,T-\tau)}L^2}\lesssim ||\nabla|^sv(t)|_{L^2}+\left(\int_t^T \frac{d\tau}{[\lambda(\tau)]^{2+\gamma'(s-\sigma_c)}}\right)^{\frac{1}{\gamma'}},$$ and the continuity $u\in \mathcal C([0,T),\dot{H}^s)$ now implies that $u(t)$ is Cauchy sequence in $\dot{H}^s$ as $t\to T$, and \fref{convustarhsigma} follows.

\begin{remark}
\label{rkintro}
 Note that the case $s=\sigma_c$ in \fref{nvbnveoeoe} leads to the logarithmic upper bound on the critical norm \fref{cnnooeen}.
 \end{remark}

{\bf step 4} Behavior of $u^*$ on the blow up point.\\

It remains to prove \fref{nonsmoothustar} which follows by adapting the argument in \cite{MR5}.\\
Let a smooth radially symmetric cut off function $\chi(r)=1$ for $r\leq 1$ and $\chi(r)=0$ for $r\geq 2$. Fix $t\in [0,T)$ and let 
\be
\label{defr}
R(t)=A_0\lambda(t)
\ee with $A_0$ given by 
\be
\label{defafao}
A_0=e^{2a\frac{\theta(2)}{b_0}}.
\ee 
We then compute the flux of $L^2$ norm: 
\bee
\left|\frac{d}{d\tau}\int\chi\left(\frac{x-x(T)}{R(t)}\right)|u(\tau)|^2\right| & = & \left|\frac{2}{R(t)}Im\left(\int \nabla \chi \left(\frac{x-x(T)}{R(t)}\right)\cdot\nabla u(\tau)\overline{u(\tau)}\right)\right|\\
& \lesssim & \frac{1}{R(t)}|u(\tau)|_{\dot{H}^{\frac{1}{2}}}^2\lesssim \frac{1}{R(t)}\frac{1}{[\lambda(\tau)]^{1-2\sigma_c}}
\eee
where we used \fref{estebootfinal}, \fref{esteisigmafinal}. We integrate this from $t$ to $T$, divide by $R^{2\sigma_c}(t)$ and get from \fref{convustarhsigma}:
\bea
\label{premooeoeo}
\nonumber & & \left|\frac{1}{R^{2\sigma_c}(t)}\int \chi\left(\frac{x-x(T)}{R(t)}\right)|u^*|^2-\frac{1}{R^{2\sigma_c}(t)}\int \chi\left(\frac{x-x(T)}{R(t)}\right)|u(t)|^2\right| \\
\nonumber & \lesssim & \frac{1}{R^{2\sigma_c+1}(t)}\int_t^T\frac{d\tau}{[\lambda(\tau)]^{1-2\sigma_c}}\lesssim  \frac{1}{A_0^{2\sigma_c+1}}\frac{1}{[\lambda(t)]^{2\sigma_c+1}(t)}\int_t^T\frac{d\tau}{[\lambda(\tau)]^{1-2\sigma_c}}\\
& \lesssim & \frac{1}{A_0^{2\sigma_c+1}b_0}\leq \frac{1}{A_0^{1/2+2\sigma_c}}
\eea
where we used the self similar speed \fref{selfsmilar} and \fref{defafao}. On the other hand, we have from \fref{decompou}:
\bee
\frac{1}{R^{2\sigma_c}(t)}\int \chi\left(\frac{x-x(T)}{R(t)}\right)|u(t)|^2 & = & \frac{\lambda^{2\sigma_c}(t)}{R^{2\sigma_c}(t)}\int \chi\left[\frac{\lambda(t)}{R(t)}\left(y+\frac{x(t)-x(T)}{\lambda(t)}\right)\right]|Q_b+\e|^2(y)dy\\
& = & \frac{1}{A_0^{2\sigma_c}}\int \chi\left[\frac{1}{A_0}\left(y+\frac{x(t)-x(T)}{\lambda(t)}\right)\right]|Q_b+\e|^2(y)dy.
\eee
Now observe from the Hardy type bound \fref{estintermjfjls}, \fref{defafao} and the bound \fref{estebootfinal} that: $$\int_{|y|\leq 10A_0}|\e|^2\lesssim A_0^3\left(\int |\nabla \e|^2+\normeldloc\right)\lesssim \Gamma_{b_0}^{\frac{1}{4}}$$ and hence \fref{covergencerate} ensures: $$
\int \chi\left[\frac{1}{A_0}\left(y+\frac{x(t)-x(T)}{\lambda(t)}\right)\right]|Q_b+\e|^2=\int |Q_p|^2(1+\delta(p))$$ with $\delta(p)\to 0$ as $p\to p_c$. Injecting this into \fref{premooeoeo} yields: 
$$\frac{1}{R^{2\sigma_c}(t)}\int \chi\left(\frac{x-x(T)}{R(t)}\right)|u^*|^2=\frac{1}{A_0^{2\sigma_c}}\int |Q_p|^2(1+\delta(p))+O\left(\frac{1}{A_0^{1/2+2\sigma_c}}\right).$$ We now let $t\to T$ ie $R(t)\to 0$ from \fref{defr} and \fref{nonsmoothustar} follows.\\
This concludes the proof of Theorem \ref{thmmain}.

\section*{Appendix}

This appendix is dedicated to the proof of \fref{cnionofiorg}. The authors are grateful to Fabrice Planchon for having put them on the right tracks.\\
After rescaling, we have:
\be
\label{correc2}
||\nabla|^{\sigma} f(\ut)|_{L^{r'}} =\frac{1}{\lambda^{p(\tilde{\sigma}-\sigma_c)}}||\nabla|^{\sigma}(F(Q_b+\e)-F(\e))|_{L^{r'}},
\ee
where the function $F:\C\goto\C$ is defined by: 
\be\label{correc3}
F(z)=|z|^{p-1}z.
\ee
Thus, \eqref{cnionofiorg} is equivalent to:
\be\label{correc4}
||\nabla|^{\sigma}(F(Q_b+\e)-F(\e))|_{L^{r'}} \lesssim ||\nabla|^{\tilde{\sigma}}\e|_{L^2}.
\ee

We now concentrate on proving \eqref{correc4}. We first rewrite $F(Q_b+\e)-F(\e)$ as:
\be\label{correc5}
F(Q_b+\e)-F(\e)=\left(\int_0^1\partial_zF(Q_b+\tau\e)d\tau\right)\e+\left(\int_0^1\partial_{\bar{z}}F(Q_b+\tau\e)d\tau\right)\bar{\e}.
\ee
Both terms in the right-hand side of \eqref{correc5} are treated in the same way. Thus, for simplicity we  may only consider the first term in the right-hand side of \eqref{correc5}. We introduce the real number $q$ such that 
\be\label{correc6}
\frac{1}{q}=\frac{1}{r'}-\frac{1}{r}.
\ee
In particular, using the definition of $r$ \eqref{strichartzparit} and usual Sobolev embeddings, we have for any function $h$:
\be\label{correc7}
|h|_{L^{(p-1)q}}\lesssim ||\nabla|^\sigma h|_{L^r}.
\ee 
Also, using again Sobolev embeddings together with the definition of $\tilde{\sigma}$ implies:
\be\label{correc8}
||\nabla|^\sigma h|_{L^r}\lesssim ||\nabla|^{\tilde{\sigma}} h|_{L^2}.
\ee 
We also introduce the real number $\nu>1$ defined by:
\be\label{correc9}
\frac{1}{\nu}=\frac{1}{r'}-\frac{1}{(p-1)q}.
\ee
\eqref{correc6} and \eqref{correc9} together with standard commutator estimates -see \cite{KatoPonce}- yield:
\be\label{correc10}
\begin{array}{r}
\ds\left||\nabla|^{\sigma}\left(\e\int_0^1\partial_zF(Q_b+\tau\e)d\tau\right)\right|_{L^{r'}}\lesssim  \ds||\nabla|^\sigma\e|_{L^r}\left|\int_0^1\partial_zF(Q_b+\tau\e)d\tau\right|_{L^{q}}\\
 \ds +|\e|_{L^{(p-1)q}}\left||\nabla|^\sigma\left[\int_0^1\partial_zF(Q_b+\tau\e)d\tau\right]\right|_{L^{\nu}},
\end{array}
\ee
which together with \eqref{correc7} and \eqref{correc8} yields:
\be\label{correc11}
\begin{array}{r}
\ds\left||\nabla|^{\sigma}\left(\e\int_0^1\partial_zF(Q_b+\tau\e)d\tau\right)\right|_{L^{r'}}\lesssim  \ds||\nabla|^{\tilde{\sigma}}\e|_{L^2}\bigg(\int_0^1\left|\partial_zF(Q_b+\tau\e)\right|_{L^{q}}d\tau\\
 \ds +\int_0^1\left||\nabla|^\sigma\left[\partial_zF(Q_b+\tau\e)\right]\right|_{L^{\nu}}d\tau\bigg).
\end{array}
\ee
Thus, in view of \eqref{correc5} and \eqref{correc11}, proving \eqref{correc4} is equivalent to proving the following bound:
\be\label{correc12}
\ds\int_0^1\left|\partial_zF(Q_b+\tau\e)\right|_{L^{q}}d\tau+\int_0^1\left||\nabla|^\sigma\left[\partial_zF(Q_b+\tau\e)\right]\right|_{L^{\nu}}d\tau\lesssim 1.
\ee
Now, we have by homogeneity:
$$\forall \tau\in[0,1], \ \ \left|\partial_zF(Q_b+\tau\e)\right|\lesssim |Q_b|^{p-1}+|\e|^{p-1}$$
which together with \eqref{corec1}, \eqref{correc7} and \eqref{correc8} yields:
\be\label{correc13}
\begin{array}{ll}
\ds\int_0^1\left|\partial_zF(Q_b+\tau\e)\right|_{L^{q}}d\tau&\ds \lesssim \int_0^1(|Q_b|^{p-1}_{L^{q(p-1)}}+|\e|^{p-1}_{L^{q(p-1)}})d\tau\\
& \ds\lesssim \int_0^1(1+||\nabla|^{\tilde{\sigma}}\e|^{p-1}_{L^2})d\tau\\
&\ds \lesssim 1.
\end{array}
\ee
Thus, we have reduced the proof of \eqref{cnionofiorg} to the proof of the following bound:
\be\label{correc14}
\int_0^1\left||\nabla|^\sigma\left[\partial_zF(Q_b+\tau\e)\right]\right|_{L^{\nu}}d\tau\lesssim 1,
\ee
where $\nu$ is defined in \eqref{correc9}.

From now on, we concentrate on proving \eqref{correc14}. To ease the notations, we define:
\be\label{correc15}
h_\tau=Q_b+\tau\e,\,0\leq\tau\leq 1.
\ee
Recall -see \cite{Cazenavebook}- the equivalent expression of homogeneous Besov norms: $\forall 0<\tilde{\sigma}<1$, $\forall 1<q<+\infty$,
\be
\label{esthomegbesov}
 |u|_{\dot{B}^{\tilde{\sigma}}_{q,2}}\sim \left(\int _0^{+\infty}\left[t^{-\tilde{\sigma}}\sup_{|y|\leq t}|u(\cdot-y)-u(\cdot)|_{L^q}\right]^2\frac{dt}{t}\right)^{\frac{1}{2}}.
\ee
Also, recall that $||\nabla|^{\tilde{\sigma}} u|_{L^q}\lesssim |u|_{\dot{B}^{\tilde{\sigma}}_{q,2}}$. 

We start by proving \eqref{correc14} in the case $1\leq N\leq 4$. In this case, $p>2$ since $p>p_c$. 
Using the homogeneity estimate: $$\forall u,v, \ \ |\partial_zF(u)-\partial_zF(v)|\lesssim |u-v|(|u|^{p-2}+|v|^{p-2}),$$ and the relation from \eqref{correc6}-\eqref{correc9}:
$$\frac{1}{\nu}=\frac{1}{r'}-\frac{1}{(p-1)q}=\frac{1}{r}+\frac{p-2}{q(p-1)},$$
we first estimate from H\"older and \fref{correc7}, \fref{correc8}:
\bee
& & |\partial_zF(h_\tau)(\cdot-y)-\partial_zF(h_\tau)(\cdot)|_{L^{\nu}} \\
& \lesssim & \left|(h_\tau(\cdot-y)-h_\tau(\cdot))[|h_\tau(\cdot-y)|^{p-2}+|h_\tau(\cdot)|^{p-2}]\right|_{L^{\nu}}\\
& \lesssim & \left|h_\tau(\cdot-y)-h_\tau(\cdot)\right|_{L^{r}}|h_\tau|^{p-2}_{L^{q(p-1)}}\\
& \lesssim & \left|h_\tau(\cdot-y)-h_\tau(\cdot)\right|_{L^{r}}||\nabla|^{\tilde{\sigma}} h_\tau|^{p-2}_{L^2}
\eee
and hence from \eqref{correc8} and \fref{esthomegbesov}: 
\be\label{correc16}
||\nabla|^\sigma\left[\partial_zF(h_\tau)\right]|_{L^{\nu}}\lesssim ||\nabla|^{\sigma} h_\tau|_{L^r}||\nabla|^{\tilde{\sigma}} h_\tau|^{p-2}_{L^2}\lesssim ||\nabla|^{\tilde{\sigma}} h_\tau|^{p-1}_{L^2}.
\ee
\eqref{corec1} and the definition of $h_\tau$ \eqref{correc15} yield:
\be\label{correc17}
||\nabla|^{\tilde{\sigma}} h_\tau|_{L^2}\lesssim ||\nabla|^{\tilde{\sigma}} Q_b|_{L^2}+||\nabla|^{\tilde{\sigma}}\e|_{L^2}\lesssim 1,
\ee
which together with \eqref{correc16} implies the wanted estimate \eqref{correc14}.

We turn to the proof of \eqref{correc14} in the remaining case $N=5$. In this case, $p>9/5$. We define the real number $\theta$ by:
\be\label{correc18}
\frac{1}{\theta}=\frac{1}{\nu}-\frac{p-9/5}{q(p-1)}.
\ee
Using the homogeneity estimate: $$\forall u,v, \ \ |\partial_zF(u)-\partial_zF(v)|\lesssim |u-v|^{\frac{4}{5}}(|u|^{p-\frac{9}{5}}+|v|^{p-\frac{9}{5}}),$$ we first estimate from H\"older, \eqref{correc7}, \eqref{correc8} and \eqref{correc17}:
\be\label{correc19}
\begin{array}{lll}
& & |\partial_zF(h_\tau)(\cdot-y)-\partial_zF(h_\tau)(\cdot)|_{L^{\nu}} \\
& \lesssim & \left|(h_\tau(\cdot-y)-h_\tau(\cdot))^{\frac{4}{5}}[|h_\tau(\cdot-y)|^{p-\frac{9}{5}}+|h_\tau(\cdot)|^{p-\frac{9}{5}}]\right|_{L^{\nu}}\\
& \lesssim & \left|h_\tau(\cdot-y)-h_\tau(\cdot)\right|^{\frac{4}{5}}_{L^{\frac{4\theta}{5}}}|h_\tau|^{p-\frac{9}{5}}_{L^{q(p-1)}}\\
& \lesssim & \left|h_\tau(\cdot-y)-h_\tau(\cdot)\right|^{\frac{4}{5}}_{L^{\frac{4\theta}{5}}}||\nabla|^{\tilde{\sigma}}h_\tau|^{p-\frac{9}{5}}_{L^2}\\
& \lesssim & \left|h_\tau(\cdot-y)-h_\tau(\cdot)\right|^{\frac{4}{5}}_{L^{\frac{4\theta}{5}}}.
\end{array}
\ee
We decompose the integral in \fref{esthomegbesov} in $t\geq 1$ and $t\leq 1$. For $t\geq 1$, we use:
$$\left|h_\tau(\cdot-y)-h_\tau(\cdot)\right|^{\frac{4}{5}}_{L^{\frac{4\theta}{5}}}\leq 1+\left|h_\tau(\cdot-y)-h_\tau(\cdot)\right|_{L^{\frac{4\theta}{5}}}$$
so that:
\be\label{correc20}
\begin{array}{l}
\ds\left(\int _1^{+\infty}\left[t^{-\sigma}\sup_{|y|\leq t}|h_\tau(\cdot-y)-h_\tau(\cdot)|^{\frac{4}{5}}_{L^{\frac{4\theta}{5}}}\right]^2\frac{dt}{t}\right)^{\frac{1}{2}}\\
\ds\lesssim \left(\int _1^{+\infty}\frac{dt}{t^{1+2\sigma}}\right)^{\frac{1}{2}}+|h_\tau|_{\dot{B}^{\sigma}_{\frac{4\theta}{5},2}}\ds\lesssim 1+|h_\tau|_{\dot{B}^{\sigma}_{\frac{4\theta}{5},2}}.
\end{array}
\ee
By usual Sobolev embeddings, we have:
\be\label{correc21}
|h_\tau|_{\dot{B}^{\sigma}_{\frac{4\theta}{5},2}}\lesssim ||\nabla|^{\sigma_1}h_\tau|_{L^2},
\ee
where $\sigma_1$ is defined by 
$$\frac{1}{2}-\frac{\sigma_1}{5}=\frac{5}{4\theta}-\frac{\sigma}{5}.$$
Using the definition of $r$, $\gamma$, $q$, $\nu$, $\theta$ and $\sigma_1$, we obtain:
$$\sigma_1=\frac{2}{\gamma}+\frac{3\sigma}{4}$$
which satisfies $\sigma<\sigma_1<1$ for $p$ close enough to $p_c$ and $\sigma$ small enough. Thus,
 the bootstrap assumptions \eqref{estenboot} \eqref{esteisigmaboot} and the definition of $h$ \eqref{correc15} yield:
\be\label{correc22}
||\nabla|^{\sigma_1} h_\tau|_{L^2}\lesssim ||\nabla|^{\sigma_1}Q_b|_{L^2}+||\nabla|^{\sigma_1}\e|_{L^2}\lesssim 1,
\ee
which together with \eqref{correc20} and \eqref{correc21} implies:
\be\label{correc23}
\begin{array}{l}
\ds\left(\int _1^{+\infty}\left[t^{-\sigma}\sup_{|y|\leq t}|h_\tau(\cdot-y)-h_\tau(\cdot)|^{\frac{4}{5}}_{L^{\frac{4\theta}{5}}}\right]^2\frac{dt}{t}\right)^{\frac{1}{2}}\lesssim 1.
\end{array}
\ee

For $t\leq 1$, we use:
$$\left|h_\tau(\cdot-y)-h_\tau(\cdot)\right|^{\frac{4}{5}}_{L^{\frac{4\theta}{5}}}\leq t^{5\sigma}+t^{-\frac{5\sigma}{4}}\left|h_\tau(\cdot-y)-h_\tau(\cdot)\right|_{L^{\frac{4\theta}{5}}},$$
so that:
\be\label{correc24}
\begin{array}{l}
\ds\left(\int _0^{1}\left[t^{-\sigma}\sup_{|y|\leq t}|h_\tau(\cdot-y)-h_\tau(\cdot)|^{\frac{4}{5}}_{L^{\frac{4\theta}{5}}}\right]^2\frac{dt}{t}\right)^{\frac{1}{2}}\\
\ds\lesssim \left(\int _0^{1}\frac{dt}{t^{1-8\sigma}}\right)^{\frac{1}{2}}+|h_\tau|_{\dot{B}^{\frac{9\sigma}{4}}_{\frac{4\theta}{5},2}}\ds\lesssim 1+|h_\tau|_{\dot{B}^{\frac{9\sigma}{4}}_{\frac{4\theta}{5},2}}.
\end{array}
\ee
By usual Sobolev embeddings, we have:
\be\label{correc25}
|h_\tau|_{\dot{B}^{\frac{9\sigma}{4}}_{\frac{4\theta}{5},2}}\lesssim ||\nabla|^{\sigma_2}h_\tau|_{L^2},
\ee
where $\sigma_2$ is defined by 
$$\frac{1}{2}-\frac{\sigma_2}{5}=\frac{5}{4\theta}-\frac{9\sigma}{20}.$$
Using the definition of $r$, $\gamma$, $q$, $\nu$, $\theta$ and $\sigma_2$, we obtain:
$$\sigma_2=\frac{2}{\gamma}+2\sigma$$
which satisfies $\sigma<\sigma_2<1$ for $p$ close enough to $p_c$ and $\sigma$ small enough. Thus,
 the bootstrap assumptions \eqref{estenboot} \eqref{esteisigmaboot} and the definition of $h$ \eqref{correc15} yield:
\be\label{correc26}
||\nabla|^{\sigma_2} h_\tau|_{L^2}\lesssim ||\nabla|^{\sigma_2}Q_b|_{L^2}+||\nabla|^{\sigma_2}\e|_{L^2}\lesssim 1,
\ee
which together with \eqref{correc24} and \eqref{correc25} implies:
\be\label{correc27}
\begin{array}{l}
\ds\left(\int _0^{1}\left[t^{-\sigma}\sup_{|y|\leq t}|h_\tau(\cdot-y)-h_\tau(\cdot)|^{\frac{4}{5}}_{L^{\frac{4\theta}{5}}}\right]^2\frac{dt}{t}\right)^{\frac{1}{2}}\lesssim 1.
\end{array}
\ee
Finally, \eqref{esthomegbesov}, \eqref{correc19}, \eqref{correc23} and \eqref{correc27} yield:
\be\label{correc28}
||\nabla|^\sigma\left[\partial_zF(h_\tau)\right]|_{L^{\nu}}\lesssim 1,
\ee
which implies the wanted estimate \eqref{correc14}. This concludes the proof of \fref{cnionofiorg}.\\


\begin{thebibliography}{10}

 %
 \bibitem{Cazenavebook} Cazenave, T.; Semilinear Schr\"odinger equations. Courant Lecture Notes in Mathematics, 10. New York University, Courant Institute of Mathematical Sciences, New York; American Mathematical Society, Providence, RI, 2003.
 
%
\bibitem{spectrenls} Chang, S.M.; Gustafson, S.; Nakanishi, K.; Tsai, T-P., Spectra of linearized operators for NLS solitary waves. SIAM J. Math. Anal. 39 (2007/08), no. 4, 1070--1111.

%
\bibitem{FGW} Fibich, G.; Gavish, N.; Wang, X.P., Singular ring solutions of critical and supercritical nonlinear Schr\"odinger equations, Physica D: Nonlinear Phenomena, 231 (2007), no. 1, 55--86.

%
\bibitem{FMR} Fibich, G.; Merle, F.; Raphael, P., Numerical proof of a spectral property related to the singularity formation for the $L^2$ critical nonlinear Schr\"odinger equation,  Phys. D  220  (2006),  no. 1, 1--13.

%
\bibitem{FP} Fibich, G.; Papanicolaou, G., Self-focusing in the perturbed and unperturbed nonlinear Schr\"odinger equation in critical dimension, SIAM J. Appl. Math. 60 (2000), no. 1, 183--240.

%
\bibitem{GN} Gidas, B.; Ni, W.M.; Nirenberg, L.,
Symmetry and related properties via the maximum principle,
Comm. Math. Phys. {\bf 68} (1979), 209---243.

%
\bibitem{GV} Ginibre, J.; Velo, G., On a class of nonlinear Schr\"odinger equations. I. The Cauchy problem, general case, J. Funct. Anal. 32 (1979), no. 1, 1--32. 

%
\bibitem{KatoPonce} Kato, T.; Ponce, G., Commutator estimates and the Euler and Navier-Stokes equations, Comm. Pure Appl. Math. 41 (1988), no. 7, 891--907.


%
\bibitem{landman} Kopell, N.; Landman, M., Spatial structure of the focusing singularity of the nonlinear Schr\"odinger equation: a geometrical analysis, SIAM J. Appl. Math. 55 (1995), no. 5, 1297--1323. 

%
\bibitem{KMR} Krieger, J.; Martel, Y.; Raphael, P., Two soliton solutions to the gravitational Hartree equation, to appear in Comm. Pure and App. Math.

%
\bibitem{KW} Kwong, M. K., Uniqueness of positive solutions of $\Delta u-u+u\sp p=0$ in ${R}\sp n$. Arch. Rational Mech. Anal. 105 (1989), no. 3, 243--266.%

%
\bibitem{MR1} Merle, F.; Rapha\"el, P., Blow up dynamic and upper bound on the blow up rate for critical nonlinear Schr\"odinger equation, Ann. Math. 161 (2005), no. 1, 157--222.

%
\bibitem{MR2} Merle, F.; Rapha\"el, P., Sharp upper bound on the blow up rate for critical nonlinear Schr\"odinger equation, Geom. Funct. Anal. 13 (2003), 591-642.

%
\bibitem{MR3} Merle, F.; Rapha\"el, P., On universality of blow up profile for $L^2$ critical nonlinear Schr\"odinger equation, Invent. Math. 156, 565-672 (2004).

%
\bibitem{MR4} Merle, F.; Rapha\"el, P., Sharp lower bound on the blow up rate for critical nonlinear Schr\"odinger equation, J. Amer. Math. Soc. 19 (2006), no. 1, 37--90.

%
\bibitem{MR5} Merle, F.; Rapha\"el, P.,  Profiles and quantization of the blow up mass for critical nonlinear Schr\"odinger equation, Comm. Math. Phys.  253  (2005),  no. 3, 675--704.

%
\bibitem{MR7} Merle, F.; Rapha\"el, Pierre, Blow up of the critical norm for some radial $L\sp 2$ super critical nonlinear Schr\"odinger equations, Amer. J. Math. 130 (2008), no. 4, 945--978.


%
\bibitem{P} Perelman, G., On the blow up phenomenon for the critical nonlinear Schr\"odinger equation in 1D, Ann. Henri. Poincaré, 2 (2001), 605-673.

%
\bibitem{R1} Rapha\"el, P., Stability of the log-log bound for blow up solutions to the critical nonlinear Schr\"odinger equation, Math. Ann. 331 (2005), 577--609.

%
\bibitem{R2} Rapha\"el, P., Existence and stability of a solution blowing up on a sphere for a $L^2$ supercritical nonlinear Schr\"odinger equation, Duke Math. J. 134 (2006), no. 2, 199--258.

%
\bibitem{RS} Rapha\"el, P., Szeftel, J., Standing ring blow up solutions to the quintic NLS in dimension $N$, to appear in Comm. Math. Phys.

%
\bibitem{RRS} Rapha\"el, P.; Rodnianski, I., Stable blow up dynamics for the critical Wave Maps and Yang-Mills, in preparation.

%
\bibitem{RodSter} Rodnianski, I.; Sterbenz, J., On the singularity formation for the critical $O(3)$ $\sigma$ model, to appear in Annals of Math.



%
\bibitem{SS} Sulem, C.; Sulem, P.L., The nonlinear Schr\"odinger equation. Self-focusing and wave collapse. Applied Mathematical Sciences, 139. Springer-Verlag, New York,
1999.

\bibitem{weinstein}  Weinstein, M.I., Modulational stability of ground 
states of nonlinear Schr\"odinger equations, SIAM J. Math. Anal.
{\bf 16} (1985), 472---491.

%
\bibitem{ZS} Zakharov, V.E.; Shabat, A.B., Exact theory
of two-dimensional self-focusing and one-dimensional self-modulation
of waves in non-linear media, Sov. Phys. JETP 34 (1972),
62--69.

\end{thebibliography}
\end{document}